# ESTIMATION OF A SEMIPARAMETRIC
# TRANSFORMATION MODEL

By Oliver Linton,[1] Stefan Sperlich[2] and Ingrid Van Keilegom[3]

*London School of Economics, Georg-August Universität Göttingen and
Université catholique de Louvain*

This paper proposes consistent estimators for transformation parameters in semiparametric models. The problem is to find the optimal transformation into the space of models with a predetermined regression structure like additive or multiplicative separability. We give results for the estimation of the transformation when the rest of the model is estimated non- or semi-parametrically and fulfills some consistency conditions. We propose two methods for the estimation of the transformation parameter: maximizing a profile likelihood function or minimizing the mean squared distance from independence. First the problem of identification of such models is discussed. We then state asymptotic results for a general class of nonparametric estimators. Finally, we give some particular examples of nonparametric estimators of transformed separable models. The small sample performance is studied in several simulations.

**1. Introduction.** Taking transformations of the data has been an integral part of statistical practice for many years. Transformations have been used to aid interpretability as well as to improve statistical performance. An important contribution to this methodology was made by Box and Cox (1964) who proposed a parametric power family of transformations that nested the logarithm and the level. They suggested that the power transformation, when applied to the dependent variable in a linear regression setting, might induce normality, error variance homogeneity and additivity

Received May 2006; revised April 2007.
[1]Supported by the ESRC.
[2]Supported by the Spanish DGI of the ministry CyT Grant SEJ2004-04583/ECON.
[3]Supported by IAP research networks nr. P5/24 and P6/03 of the Belgian government (Belgian Science Policy).

of effects. They proposed estimation methods for the regression and transformation parameters. Carroll and Ruppert (1984) applied this and other transformations to both dependent and independent variables. A number of other dependent variable transformations have been suggested, for example, the Zellner–Revankar (1969) transform and the Bickel and Doksum (1981) transform. The transformation methodology has been quite successful and a large literature exists on this subject for parametric models; see Carroll and Ruppert (1988). In survival analysis there are many applications due to the interpretation of versions of the model as accelerated failure time models, proportional hazard models, mixed proportional hazard models and proportional odds models; see, for example, Doksum (1987), Wei (1992), Cheng and Wu (1994), Cheng, Wei and Ying (1995) and van den Berg (2001).

In this work we concentrate on transformations in a regression setting. For many data, linearity of covariate effect after transformation may be too strong. We consider a rather general specification, allowing for nonparametric covariate effects. Let $X$ be a $d$-dimensional random vector and $Y$ be a random variable, and let $\{(X_i, Y_i)\}_{i=1}^n$ be an i.i.d. sample from this population. Consider the estimation of the regression function $m(x) = E(Y \mid X = x)$. Stone (1980, 1982) and Ibragimov and Hasminskii (1980) showed that the optimal rate for estimating $m$ is $n^{-\ell/(2\ell+d)}$, with $\ell$ a measure of the smoothness of $m$. This rate of convergence can be very slow for large dimensions $d$. One way of achieving better rates of convergence is making use of dimension reducing separability structures. The most common examples are additive or multiplicative modeling. An additive structure for $m$, for example, is a regression function of the form $m(x) = \sum_{\alpha=1}^d m_\alpha(x_\alpha)$, where $x = (x_1, \ldots, x_d)^\top$ are the $d$-dimensional predictor variables and $m_\alpha$ are one-dimensional nonparametric functions. Stone (1986) showed that for such regression curves the optimal rate for estimating $m$ is the one-dimensional rate of convergence $n^{-\ell/(2\ell+1)}$. Thus, one speaks of dimensionality reduction through additive modeling.

We examine a semiparametric model that combines a parametric transformation with the flexibility of an additive nonparametric regression function. Suppose that

$$(1) \qquad \Lambda(Y) = G(m_1(X_1), \ldots, m_d(X_d)) + \varepsilon,$$

where $\varepsilon$ is independent of $X$, while $G$ is a known function and $\Lambda$ is a monotonic function. Special cases of $G$ are $G(z) = H(\sum_{\alpha=1}^d z_\alpha)$ and $G(z) = H(\prod_{\alpha=1}^d z_\alpha)$ for some strictly monotonic known function $H$. The general model in which $\Lambda$ is monotonic and $G(z) = \sum_{\alpha=1}^d z_\alpha$ was previously addressed in Breiman and Friedman (1985) who suggested estimation procedures based on the iterative backfitting method, which they called ACE. However, they did not provide many results about the statistical properties of their procedures. Linton, Chen, Wang and Härdle (1997) considered the model with $\Lambda = \Lambda_\theta$ parametric and additive $G$, $G(z) = \sum_{\alpha=1}^{d} z_\alpha$. They proposed to estimate the parameters of the transformation $\Lambda$ by either an instrumental variable method or a pseudo-likelihood method based on Gaussian $\varepsilon$. For the instrumental variable method, they assumed that identification held from some unconditional moment restriction but they did not provide justification for this from primitive conditions. Unfortunately, our simulation evidence suggests that both methods work poorly in practice and may even be inconsistent for many parameter configurations. To estimate the unknown functions $m_\alpha$ they used the marginal integration method of Linton and Nielsen (1995) and, consequently, their method cannot achieve the semiparametric efficiency bound for estimation of $\theta$ even in the few cases where Gaussian errors are well defined and their method is consistent.

We argue that an even more general version of the model (1) is identified following results of Ekeland, Heckman and Nesheim (2004). For practical reasons, we propose estimation procedures only for the parametric transformation case where $\Lambda(y) = \Lambda_{\theta_o}(y)$ for some parametric family $\{\Lambda_\theta(\cdot), \theta \in \Theta\}$ of transformations where $\Theta \subset \mathbb{R}^k$. This model includes, for example, the Nielsen, Linton and Bickel (1998) (reversed) proportional hazard model where the baseline hazard is parametric and the covariate effect is nonparametric. This is appropriate for certain mortality studies where there are well established models for baseline mortality but covariate effects are not so well understood. To estimate the transformation parameters, we use two approaches. First, a semiparametric profile likelihood estimator (PL) that involves nonparametric estimation of the density of $\varepsilon$, and second, a mean squared distance from the independence method (MD) based on estimated c.d.f.'s of $(X, \varepsilon)$. Both methods use a profiled estimate of the (separable) nonparametric components of $m_\theta$. We use both the integration method and the smooth backfitting method of Mammen, Linton and Nielsen (1999) to estimate these components. The MD estimator involves discontinuous functions of nonparametric estimators and we use the theory of Chen, Linton and Van Keilegom (2003) to obtain its asymptotic properties. We derive the asymptotic distributions of our estimators under standard regularity conditions, and we show that the estimators of $\theta_o$ are root-$n$ consistent. The corresponding estimators of the component functions $m_j(\cdot)$ behave as if the parameters $\theta_o$ were known and are also asymptotically normal at nonparametric rates.

The rest of the paper is organized as follows. In the next section we clarify identification issues. In Section 3 we introduce the two estimators for the transformation parameter. Section 4 contains the asymptotic theory of these two estimators. Additionally, we discuss tools like bootstrap for possible inference on the transformation parameter. Finally, in Section 5 we study the finite sample performance of all methods presented and compare the different estimators of the transformation parameter, as well as the different estimators of the additive components in this context. A special emphasis is also given to the question of bandwidth choice. All proofs are deferred to Appendix A and Appendix B.

## 2. Nonparametric identification.   Suppose that

$$\Lambda(Y) = m(X) + \varepsilon, \tag{2}$$

where $\varepsilon$ is independent of $X$ with unknown distribution $F_\varepsilon$, and the functions $\Lambda$ and $m$ are unknown. Then

$$F_{Y|X}(y,x) = \Pr[Y \le y | X = x] = F_\varepsilon(\Lambda(y) - m(x)). \tag{3}$$

Recently Ekeland, Heckman and Nesheim (2004), building on ideas of Horowitz (1996, 2001), have shown that this model is identifiable up to a couple of normalizations under smoothness conditions on $(F_\varepsilon, \Lambda, m)$ and monotonicity conditions on $\Lambda$ and $F_\varepsilon$. The basic idea is to note that, for each $j$,

$$\frac{\partial F_{Y|X}(y,x)}{\partial y} \Big/ \frac{\partial F_{Y|X}(y,x)}{\partial x_j} = -\frac{\lambda(y)}{\partial m(x)/\partial x_j}, \tag{4}$$

where $\lambda(y) = \partial \Lambda(y)/\partial y$. Then by integrating out either $y$ or $x$, one obtains $\lambda(\cdot)$ up to a constant or $\partial m(\cdot)/\partial x_j$ up to a constant. By further integrations, one obtains $\Lambda(\cdot)$ and $m(\cdot)$ up to a constant. One then obtains $F_\varepsilon$ by inverting the relationship (3) and imposing the normalizations. Horowitz (1996) indeed covers the special case where $m(x)$ is linear.

The above arguments show that for identification it is not necessary to restrict $\Lambda$, $m$ or $F_\varepsilon$ beyond monotonicity, smoothness and normalization restrictions. However, the implied estimation strategy can be very complicated; see, for example, Lewbel and Linton (2006). In addition, the fully nonparametric model does not at all reduce the curse of dimensionality in comparison with the unrestricted conditional distribution $F_{Y|X}(y,x)$, which makes the practical relevance of the identification result limited. This is why we consider additive and multiplicative structures on $m$ and a parametric restriction on $\Lambda$. The unrestricted model could be used for testing of these assumptions, although we do not pursue this in this paper.

To conclude this section, we discuss briefly some related work on identification of related models. Linton, Chen, Wang and Härdle (1997) assumed identification of the model (2) with parametric $\Lambda$ and additive $m$ based on an unconditional moment restriction on the error term rather than full independence. In particular, they assumed that $E[Z\varepsilon] = 0$ for a vector of variables $Z$. This does not seem to be sufficient to justify identification and, indeed, our simulation evidence supports this concern. Finally, we mention a nonparametric identification result of Breiman and Friedman (1985). They defined functions $\Lambda(\cdot), m_1(\cdot), \ldots, m_d(\cdot)$ as minimizers of the least squares objective function

$$(5) \qquad e^2(\Lambda, m_1, \ldots, m_d) = \frac{E[\{\Lambda(Y) - \sum_{\alpha=1}^{d} m_\alpha(X_\alpha)\}^2]}{E[\Lambda^2(Y)]}$$

for general random variables $Y, X_1, \ldots, X_d$. They showed the existence of minimizers of (5) and showed that the set of minimizers forms a finite dimensional linear subspace (of an appropriate class of functions) under additional conditions. These conditions were that: (i) $\Lambda(Y) - \sum_{\alpha=1}^{d} m_\alpha(X_\alpha) = 0$ a.s. implies that $\Lambda(Y), m_\alpha(X_\alpha) = 0$ a.s., $\alpha = 1, \ldots, d$; (ii) $E[\Lambda(Y)] = 0, E[m_\alpha(X_\alpha)] = 0$, $E[\Lambda^2(Y)] < \infty$, and $E[m_\alpha^2(X_\alpha)] < \infty$; (iii) The conditional expectation operators $E[\Lambda(Y)|X_\alpha]$, $E[m_\alpha(X_\alpha)|Y]$, $\alpha = 1, \ldots, d$ are compact. This result does not require any model assumptions like conditional moments or independent errors, but has more limited scope. We shall maintain the model assumption of independent errors in the sequel.

**3. Estimating the transformation.** In the sequel we consider the model

$$(6) \qquad \Lambda_{\theta_o}(Y) = m(X) + \varepsilon,$$

where $\{\Lambda_\theta : \theta \in \Theta\}$ is a parametric family of strictly increasing functions, while the function $m(\cdot)$ is of unknown form but with a certain predetermined structure that is sufficient to yield dimensionality reduction. We assume that the error term $\varepsilon$ is independent of $X$, has distribution $F$, and $E(\varepsilon) = 0$. The covariate $X$ is $d$-dimensional and has compact support $\mathcal{X} = \prod_{\alpha=1}^{d} R_{X_\alpha}$. Among the many transformations of interest, the following ones are used most commonly: (Box–Cox) $\Lambda_\theta(y) = \frac{y^\theta - 1}{\theta}$ ($\theta \neq 0$) and $\Lambda_\theta(y) = \log(y)$ ($\theta = 0$); (Zellner–Revankar) $\Lambda_\theta(y) = \ln y + \theta y^2$; (Arcsinh) $\Lambda_\theta(y) = \sinh^{-1}(\theta y)/\theta$. The arcsinh transform is discussed in Johnson ([1949](#)) and more recently in Robinson ([1991](#)). The main advantage of the arcsinh transform is that it works for $y$ taking any value, while the Box–Cox and the Zellner–Revankar transforms are only defined if $y$ is positive. For these transformations, the error term cannot be normally distributed except for a few isolated parameters, and so the Gaussian likelihood is misspecified. In fact, as Amemiya and Powell ([1981](#)) point out, the resulting estimators (in the parametric case) are inconsistent when only $n \to \infty$.

We let $\Theta$ denote a finite dimensional parameter set (a compact subset of $\mathbb{R}^k$) and $\mathcal{M}$ an infinite dimensional parameter set. We assume that $\mathcal{M}$ is a vector space of functions endowed with metric $\|\cdot\|_{\mathcal{M}} = \|\cdot\|_\infty$. We denote $\theta_o \in \Theta$ and $m_o \in \mathcal{M}$ as the true unknown finite and infinite dimensional parameters. Define the regression function

$$m_\theta(x) = E[\Lambda_\theta(Y)|X = x]$$

for each $\theta \in \Theta$. Note that $m_{\theta_o}(\cdot) \equiv m_o(\cdot)$.

We suppose that we have a randomly drawn sample $Z_i = (X_i, Y_i)$, $i = 1, \ldots, n$, from model (6). Define, for $\theta \in \Theta$ and $m \in \mathcal{M}$,

$$\varepsilon(\theta, m) = \Lambda_\theta(Y) - m(X),$$

and let $\varepsilon_\theta = \varepsilon(\theta) = \varepsilon(\theta, m_\theta)$ and $\varepsilon_o = \varepsilon_{\theta_o}$. When there is no ambiguity, we also use the notation $\varepsilon$ and $m$ to indicate $\varepsilon_o$ and $m_o$. Moreover, let $\Lambda_o = \Lambda_{\theta_o}$.

In the sequel we will denote by $\widehat{m}_\theta$ any estimator of $m_\theta$ under either the additive or the multiplicative model. In the simulation section we will focus on the additive model and the smooth backfitting estimator, denoted by $\widehat{m}_\theta^{\mathrm{BF}}(\cdot)$. See Mammen, Linton and Nielsen (1999) for its definition. $\widehat{m}_\theta^{\mathrm{BF}}$ consistently estimates a function $m_\theta^{\mathrm{BF}}(\cdot)$, where $m_{\theta_0}^{\mathrm{BF}}(\cdot) = m_{\theta_0}(\cdot)$, but $m_\theta^{\mathrm{BF}}(\cdot) \neq m_\theta(\cdot)$ for $\theta \neq \theta_0$.

### 3.1. *The profile likelihood (PL) estimator.*

The method of profile likelihood has already been applied to many different semiparametric estimation problems. The basic idea is simply to replace all unknown expressions of the likelihood function by their nonparametric (kernel) estimates. We consider $\Lambda_\theta(Y) = m_\theta(X) + \varepsilon_\theta$ for any $\theta \in \Theta$. Then, the cumulative distribution function is

$$\begin{aligned}
\Pr[Y \leq y | X] &= \Pr[\Lambda_\theta(Y) \leq \Lambda_\theta(y) | X] \\
&= \Pr[\varepsilon_\theta \leq \Lambda_\theta(y) - m_\theta(X) | X] \\
&= F_{\varepsilon(\theta)}(\Lambda_\theta(y) - m_\theta(X)),
\end{aligned}$$

where $F_{\varepsilon(\theta)}(e) = F_{\varepsilon(\theta, m_\theta)}(e)$ and $F_{\varepsilon(\theta, m)} = P(\varepsilon(\theta, m) \leq e)$, and so

$$f_{Y|X}(y|x) = f_{\varepsilon(\theta)}(\Lambda_\theta(y) - m_\theta(x)) \Lambda_\theta'(y),$$

where $f_{\varepsilon(\theta)}$ and $f_{Y|X}$ are the probability density functions of $\varepsilon(\theta)$ and of $Y$ given $X$. Then, the log likelihood function is

$$\sum_{i=1}^n \{ \log f_{\varepsilon(\theta)}(\Lambda_\theta(Y_i) - m_\theta(X_i)) + \log \Lambda_\theta'(Y_i) \}.$$

Let

$$(7) \qquad \widehat{f}_{\varepsilon(\theta)}(e) := \frac{1}{ng} \sum_{i=1}^n K_2 \left( \frac{e - \widehat{\varepsilon}_i(\theta)}{g} \right),$$

with $\widehat{\varepsilon}_i(\theta) = \widehat{\varepsilon}_i(\theta, m_\theta)$ and $\widehat{\varepsilon}_i(\theta, m) = \varepsilon_i(\theta, \widehat{m}) = \Lambda_\theta(Y_i) - \widehat{m}(X_i)$. Here, $K_2$ is a scalar kernel and $g$ is a bandwidth sequence. Then, define the profile likelihood estimator of $\theta_o$ by

$$(8) \qquad \widehat{\theta}_{\mathrm{PL}} = \operatorname*{arg\,max}_{\theta \in \Theta} \sum_{i=1}^n [\log \widehat{f}_{\varepsilon(\theta)}(\Lambda_\theta(Y_i) - \widehat{m}_\theta(X_i)) + \log \Lambda_\theta'(Y_i)].$$

The computation of $\widehat{\theta}_{\mathrm{PL}}$ can be done by grid search in the scalar case and using derivative-based algorithms in higher dimensions, assuming that the kernels are suitably smooth.

3.2. *Mean square distance from independence (MD) estimator.* There are four good reasons why it is worth providing alternative estimators when it comes to practical work. First, as we will see in Section 5, the profile likelihood method is computationally quite expensive. In particular, so far we have not found a reasonable implementation for the recentered bootstrap. Second, for that approach we do not only face the typical question of bandwidth choice for the nonparametric part $m_\theta$, we additionally face a bandwidth for the density estimation; see equation (7). Third, there are some transformation models $\Lambda_\theta$ for which the support of $Y$ depends on the parameter $\theta$ and so are nonregular. Finally, although the estimator we get from the profile likelihood is under certain conditions efficient in the asymptotic sense [Severini and Wong (1992)], this tells us little about its finite sample performance, neither in absolute terms nor in comparison with competitors.

One possible and computationally attractive competitor is the minimization of the mean square distance from independence. Why it is computationally more attractive will be explained in Section 5. This method we will introduce here has been reviewed in Koul (2001) for other problems.

Define, for each $\theta \in \Theta$ and $m \in \mathcal{M}$, the empirical distribution functions

$$\widehat{F}_X(x) = \frac{1}{n} \sum_{i=1}^{n} 1(X_i \leq x);$$

$$\widehat{F}_{\varepsilon(\theta)}(e) = \frac{1}{n} \sum_{i=1}^{n} 1(\widehat{\varepsilon}_i(\theta) \leq e);$$

$$\widehat{F}_{X,\varepsilon(\theta)}(x,e) = \frac{1}{n} \sum_{i=1}^{n} 1(X_i \leq x) 1(\widehat{\varepsilon}_i(\theta) \leq e),$$

the moment function

$$G_{n\mathrm{MD}}(\theta, \widehat{m}_\theta)(x,e) = \widehat{F}_{X,\varepsilon(\theta)}(x,e) - \widehat{F}_X(x) \widehat{F}_{\varepsilon(\theta)}(e)$$

and the criterion function

$$(9) \qquad \|G_{n\mathrm{MD}}(\theta, \widehat{m}_\theta)\|_2^2 = \int [G_{n\mathrm{MD}}(\theta, \widehat{m}_\theta)(x,e)]^2 \, d\mu(x,e)$$

for some probability measure $\mu$. We define an estimator of $\theta$, denoted $\widehat{\theta}_{\mathrm{MD}}$, as any approximate minimizer of $\|G_{n\mathrm{MD}}(\theta, \widehat{m}_\theta)\|_2^2$ over $\Theta$. To be precise, let

$$\|G_{n\mathrm{MD}}(\widehat{\theta}_{\mathrm{MD}}, \widehat{m}_{\widehat{\theta}})\|_2 = \inf_{\theta \in \Theta} \|G_{n\mathrm{MD}}(\theta, \widehat{m}_\theta)\|_2 + o_p(1/\sqrt{n}).$$

There are many algorithms available for computing the optimum of general nonsmooth functions, for example, the Nelder–Mead, and the more recent genetic and evolutionary algorithms.

We can use in (9) the empirical measure $d\mu_n$ of $\{X_i, \widehat{\varepsilon}_i(\theta)\}_{i=1}^n$, which results in a criterion function

$$(10) \qquad Q_n(\theta) = \frac{1}{n} \sum_{i=1}^n [G_{n\mathrm{MD}}(\theta, \widehat{m}_\theta)(X_i, \widehat{\varepsilon}_i(\theta))]^2.$$

In the sequel we will denote $m_\theta$ to indicate either the function $E[\Lambda_\theta(Y)|X = \cdot]$ or the function $m_\theta^{\mathrm{BF}}$ defined above (or the population version of any other estimator of $m_\theta$). It will be clear from the context which function it represents.

## 4. Asymptotic properties.

We now discuss the asymptotic properties of our procedures. Note that although nonparametric density estimation with non- or semiparametrically constructed variables has already been considered in Van Keilegom and Veraverbeke (2002) and in Sperlich (2005), their results cannot be applied directly to our problem. The first one treated the more complex problem of censored regression models but have no additional parameter like our $\theta$. Nevertheless, as they consider density estimation with nonparametrically estimated residuals, their results come much closer to our needs than the second paper. Neither offer results on derivative estimation. As we will see now, this we need when we translate our estimation problem into the estimation framework of Chen, Linton and Van Keilegom (2003) [CLV (2003) in the sequel].

To be able to apply the results of CLV (2003) for proving the asymptotics of the profile likelihood, we need an objective function that takes its minimum at $\theta_o$. Therefore, we introduce some notation. For any function $\varphi$, we define $\dot{\varphi} := \partial \varphi / \partial \theta$ and $\dot{\widehat{\varphi}} := \partial \widehat{\varphi} / \partial \theta$, respectively. Similarly, we define for any function $\varphi$: $\varphi'(u) := \partial \varphi(u) / \partial u$ and $\widehat{\varphi}'(u) := \partial \widehat{\varphi}(u) / \partial u$, respectively. The same holds for any combination of primes and dots.

We use the abbreviated notation $s = (m, r, f, g, h)$, $s_\theta = (m_\theta, \dot{m}_\theta, f_{\varepsilon(\theta)}, f'_{\varepsilon(\theta)}, \dot{f}_{\varepsilon(\theta)})$, $s_o = s_{\theta_o}$ and $\widehat{s}_\theta = (\widehat{m}_\theta, \dot{\widehat{m}}_\theta, \widehat{f}_{\varepsilon(\theta)}, \widehat{f}'_{\varepsilon(\theta)}, \dot{\widehat{f}}_{\varepsilon(\theta)})$. Then, define for any $s = (m, r, f, g, h)$,

$$
\begin{aligned}
(11) \quad & G_{n\mathrm{PL}}(\theta, s) \\
& = n^{-1} \sum_{i=1}^n \left\{ \frac{1}{f\{\varepsilon_i(\theta, m)\}} \right. \\
& \qquad\qquad \times [g\{\varepsilon_i(\theta, m)\}\{\dot{\Lambda}_\theta(Y_i) - r(X_i)\} + h\{\varepsilon_i(\theta, m)\}] \\
& \qquad\qquad\qquad \left. + \frac{\dot{\Lambda}_\theta(Y_i)}{\Lambda'_\theta(Y_i)} \right\},
\end{aligned}
$$

and let $G_{\mathrm{PL}}(\theta, s) = E[G_{n\mathrm{PL}}(\theta, s)]$, and $\Gamma_{1\mathrm{PL}} = \frac{\partial}{\partial \theta} G_{\mathrm{PL}}(\theta, s_\theta)\big\downarrow_{\theta = \theta_o}$.

Note that $\|G_{\mathrm{PL}}(\theta, s_\theta)\|$ and $\|G_{n\mathrm{PL}}(\theta, \widehat{s}_\theta)\|$ take their minimum at $\theta_o$ and $\widehat{\theta}_{\mathrm{PL}}$ respectively (where $\|\cdot\|$ denotes the Euclidean norm). We assume in the [Appendix](#) that the estimator of the nonparametric index obeys a certain asymptotic expansion. Note that, when the index is additively separable, typical candidates are the marginal integration estimator [Tjøstheim and Auestad ([1994](#)), Linton and Nielsen ([1995](#)) and Sperlich, Tjøstheim and Yang ([2002](#)) for additive interaction models] and the smooth backfitting [Mammen, Linton and Nielsen ([1999](#)) and Nielsen and Sperlich ([2005](#))]. Both estimators obey a certain asymptotic expansion. The proof of such expansions can be found in Lemmas 6.1 and 6.2 of Mammen and Park ([2005](#)) for backfitting and in Linton et al. ([1997](#)) for marginal integration. In consequence, we obtain expansions for $\widehat{f}_{\varepsilon(\theta)}, \widehat{f}'_{\varepsilon(\theta)}, \dot{\widehat{f}}_{\varepsilon(\theta)}$.

**Theorem 4.1.** *Under Assumptions [A.1–A.1](#) given in Appendix [A](#), we have*

$$\widehat{\theta}_{\mathrm{PL}} - \theta_o = -\Gamma_{1\mathrm{PL}}^{-1} G_{n\mathrm{PL}}(\theta_o, s_o) + o_p(n^{-1/2}),$$

$$\sqrt{n}(\widehat{\theta}_{\mathrm{PL}} - \theta_o) \Longrightarrow N(0, \Omega_{\mathrm{PL}}),$$

*where $\Omega_{\mathrm{PL}} = \Gamma_{1\mathrm{PL}}^{-1} \mathrm{Var}\{G_{1\mathrm{PL}}(\theta_o, s_o)\}(\Gamma_{1\mathrm{PL}}^T)^{-1}$.*

Note that the variance of $\widehat{\theta}_{\mathrm{PL}}$ equals the variance of the estimator of $\theta_o$ that is based on the true (unknown) values of the nuisance functions $m_o, \dot{m}_o, f_\varepsilon, f'_\varepsilon$ and $\dot{f}_\varepsilon$. For the smooth backfitting, we expect that the profile likelihood estimator is semiparametrically efficient following Severini and Wong ([1992](#)); see also Linton and Mammen ([2005](#)).

We obtain the asymptotic distribution of $\widehat{\theta}_{\mathrm{MD}}$ using a modification of Theorems 1 and 2 of CLV ([2003](#)). That result applied to the case where the norm in ([9](#)) was finite dimensional, although their Theorem 1 is true as stated with the more general norm. Regarding their Theorem 2, we need to modify only condition 2.5 to take account of the fact that $G_{n\mathrm{MD}}(\theta, m_\theta)$ is a stochastic process in $(x, e)$. Let $\lambda_\theta(y) = \dot{\Lambda}_\theta(y) = \partial \Lambda_\theta(y)/\partial \theta$ and let $\lambda_o = \lambda_{\theta_o}$. We also note that

$$\frac{\partial}{\partial \theta} E[\Lambda_\theta(Y)|X]\bigg\downarrow_{\theta = \theta_o} = \int \lambda_o(\Lambda_o^{-1}(m_o(X) + e)) f_\varepsilon(e)\, de.$$

Define the matrix

$$\Gamma_{1\mathrm{MD}}(x, e) = f_\varepsilon(e) E[(1(X \le x) - F_X(x))(\lambda_o(\Lambda_o^{-1}(m_o(X) + e)) + \dot{m}_o(X))],$$

and the i.i.d. mean zero and finite variance random variables

$$\overline{U}_i = \int [1(X_i \le x) - F_X(x)][1(\varepsilon_i \le e) - F_\varepsilon(e)]\Gamma_{1\mathrm{MD}}(x, e)\, d\mu(x, e)$$

$$+ f_X(X_i) \sum_{\alpha=1}^{d} v_{o1\alpha}(X_{\alpha i}, \varepsilon_i) \int f_\varepsilon(e)(1(X_i \le x) - F_X(x))$$

$$\times \Gamma_{1\mathrm{MD}}(x, e)\, d\mu(x, e),$$

where $v_{o1\alpha}(\cdot)$ is defined in Assumption A.8 in Appendix A.

Let $V_{1\mathrm{MD}} = E[\overline{U}_i \overline{U}_i^\top]$ and $\overline{\Gamma}_{1\mathrm{MD}} = \int \Gamma_{1\mathrm{MD}}(x, e)\Gamma_{1\mathrm{MD}}^T(x, e)\, d\mu(x, e)$.

THEOREM 4.2.  *Under Assumptions* B.1–B.8 *given in Appendix* B, *we have*

$$\widehat{\theta}_{\mathrm{MD}} - \theta_o = -\overline{\Gamma}_{1\mathrm{MD}}^{-1}\overline{U}_i + o_p(n^{-1/2}),$$

$$\sqrt{n}(\widehat{\theta}_{\mathrm{MD}} - \theta_o) \Longrightarrow N(0, \Omega_{\mathrm{MD}}),$$

*where* $\Omega_{\mathrm{MD}} = \overline{\Gamma}_{1\mathrm{MD}}^{-1} V_{1\mathrm{MD}} \overline{\Gamma}_{1\mathrm{MD}}^{-1}$.

REMARKS.  1. The properties of the resulting estimators of $m$ and its components follow from standard calculations as in Linton et al. (1997), Theorem 3: the asymptotic distributions are as if the parameters $\theta_o$ were known.

2. Bootstrap standard errors. CLV (2003) proposes and justifies the use of the ordinary bootstrap. Let $\{Z_i^*\}_{i=1}^n$ be drawn randomly with replacement from $\{Z_i\}_{i=1}^n$, and let

$$G_{n\mathrm{MD}}^*(\theta, m)(x, e) = \widehat{F}_{X\varepsilon(\theta)}^*(x, e) - \widehat{F}_X^*(x)\widehat{F}_{\varepsilon(\theta)}^*(e),$$

where $\widehat{F}_{X\varepsilon(\theta)}^*, \widehat{F}_X^*(x)$ and $\widehat{F}_{\varepsilon(\theta)}^*$ are computed from the bootstrap data. Let also $\widehat{m}_\theta^*(\cdot)$ (for each $\theta$) be the same estimator as $\widehat{m}_\theta(\cdot)$, but based on the bootstrap data. Following Hall and Horowitz [(1996), page 897], it is necessary to recenter the moment condition, at least in the overidentified case. Thus, define the bootstrap estimator $\widehat{\theta}_{\mathrm{MD}}^*$ to be any sequence that satisfies

$$
\begin{aligned}
(12) \quad & \|G_{n\mathrm{MD}}^*(\widehat{\theta}_{\mathrm{MD}}^*, \widehat{m}_{\widehat{\theta}_{\mathrm{MD}}^*}^*) - G_{n\mathrm{MD}}(\widehat{\theta}_{\mathrm{MD}}, \widehat{m}_{\widehat{\theta}_{\mathrm{MD}}})\| \\
& = \inf_{\theta \in \Theta} \|G_{n\mathrm{MD}}^*(\theta, \widehat{m}_\theta^*) - G_{n\mathrm{MD}}(\widehat{\theta}_{\mathrm{MD}}, \widehat{m}_{\widehat{\theta}_{\mathrm{MD}}})\| + o_{p^*}(n^{-1/2}),
\end{aligned}
$$

where superscript $*$ denotes a probability or moment computed under the bootstrap distribution conditional on the original data set $\{Z_i\}_{i=1}^n$. The resulting bootstrap distribution of $\sqrt{n}(\widehat{\theta}_{\mathrm{MD}}^* - \widehat{\theta}_{\mathrm{MD}})$ can be shown to be asymptotically the same as the distribution of $\sqrt{n}(\widehat{\theta}_{\mathrm{MD}} - \theta_o)$, by following the same arguments as in the proof of Theorem B in CLV (2003). Similar arguments can be applied to the PL method.

3. Estimated weights. Suppose that we have estimated weights $\mu_n(x, e)$ that satisfy $\sup_{x,e} |\mu_n(x, e) - \mu(x, e)| = o_p(1)$. Then the estimator computed with the estimated weights $\mu_n(x, e)$ has the same distribution theory as the estimator that used the limiting weights $\mu(x, e)$.

4. Note that the asymptotic distributions in Theorems 4.1 and 4.2 do not depend on the details of the estimator $\hat{m}_\theta^{BF}(x)$, only on their population interpretations through

$$(13) \qquad \frac{\partial m_\theta^{BF}}{\partial \theta}(\cdot) = \underset{m \in \mathcal{M}_{add}}{\arg\min} \int \left[ \left( \frac{\partial m_\theta}{\partial \theta}(X) - m(X) \right)^2 \right] f_X(X) \, dX,$$

where

$$\mathcal{M}_{add} = \left\{ m : m(x) = \sum_{\alpha=1}^d m_\alpha(x_\alpha) \text{ for some } m_1(\cdot), \ldots, m_d(\cdot) \right\}.$$

**5. Performance in finite samples.** We consider the following data generating process:

$$(14) \qquad \Lambda_\theta(Y) = b_0 + b_1 X_1^2 + b_2 \sin(\pi X_2) + \varepsilon \sigma_e,$$

where $\Lambda_\theta$ is the Box–Cox transformation, $X_1, X_2 \sim U[-0.5, 0.5]^2$ and $\varepsilon$ drawn from $N(0, 1)$ but restricted on $[-3, 3]$. We study three different models with $b_0 = 3.0\sigma_e + b_2$ and $b_1$, $b_2$, $\sigma_e$ as follows: for model 1, we set $b_1 = 5.0$, $b_2 = 2.0$, $\sigma_e = 1.5$; for model 2, $b_1 = 3.5$, $b_2 = 1.5$, $\sigma_e = 1.0$; and for model 3, $b_1 = 2.5$, $b_2 = 1.0$, $\sigma_e = 0.5$. Parameter $\theta_o$ is set to 0.0, 0.5 and 1.0. Note that $\Lambda_\theta(Y)$ is by construction always positive in our simulations.

We estimated $\theta$ by a grid search on $[-0.5, 1.5]$ with step length 0.0625. Our implementations for estimators of the additive index follow exactly Nielsen and Sperlich (2005) for the backfitting (BF), and Hengartner and Sperlich (2005) for the marginal integration (MI). We just show results for the BF method; results for marginal integration, further details and more results on the bootstrap can be found in Sperlich, Linton and Van Keilegom (2007). BF has been chosen as we know from Sperlich, Linton and Härdle (1999) that backfitting is more reliable when predicting the whole mean function—which matters more in our context—whereas MI has some advantages when looking at the marginal impacts. We use the local constant versions with quartic kernel $K(u) = \frac{15}{16}(1 - u^2)_+^2$ and bandwidth $h_1 = h_2 = n^{-1/5}h_0$ for a large range of $h_0$-values. For the density estimator of the predicted residuals in the PL, we use Silverman's rule of thumb bandwidth in each iteration.

5.1. *Comparing PL with MD.* We first evaluate robustness against bandwidth. Table 1 gives the means and standard deviations calculated for samples of size $n = 100$ from 500 replications for each $\theta_o$ and different bandwidth. Since the parameter set $\Theta = [-0.5, 1.5]$, the simulation results for $\theta_o = 0.0$ and 1.0 are biased toward the interior of the $\Theta$. Note further that there is also an interaction between bandwidth and $\theta$ (the estimated as well as the real one) concerning the smoothness of the model: using local constant smoothers, the estimates will have more bias for larger derivatives. On the other hand, both a smaller $\theta$ and a larger $h_0$ make the model "smoother," and vice versa. We therefore study the bandwidth choice in a separate simulation.

Table 1 gives the results for any combination of model, bandwidth and method. If the error distribution is small compared to the estimation error, then the MD is expected to do worse. Indeed, even though model 3 is the smoothest model and therefore the easiest estimation problem, for the smallest error standard deviation ($\sigma_e = 0.5$), the MD does worse. In those cases the PL estimator should perform better, and so it does. It might be surprising that $\theta$ mostly gets better estimated in model 1 than in model 2 and model 3, where the nonparametric functionals are much easier to estimate. But notice that for the quality of $\hat{\theta}$ the relation between estimation error and model error is more important. This is also true for the PL method. Nevertheless, at least for small samples, none of the estimators seems to outperform uniformly the other: so the PL has mostly smaller variance, whereas MD has mostly smaller bias. As expected, for very small samples, the results depend on the bandwidth. For this reason, and due to its importance in practice, we study this problem more in detail below. We should mention that the PL method is much more expensive to calculate than the MD.

5.2. *Bandwidth choice.* Perhaps the simplest approach conceptually would be to apply plug-in bandwidths. However, this method relies on asymptotic expressions with unknown functions and parameters that are even more complicated to estimate. Furthermore, in simulations [see Sperlich, Linton and Härdle (1999) or Mammen and Park (2005)] they turned out not to work satisfactorily. Instead, we applied the cross-validation method for smooth backfitting developed in Nielsen and Sperlich (2005) and adapted to our context.

In Table 2 we give the results for minimizing the MD over $\theta \in \Theta$ choosing $h \in \mathbb{R}^d$ by cross validation. Notice that we allow for different bandwidths for each additive component. The simulations are done as before, but only for model 1 and based on just 100 simulation runs what is enough to see the following: The results presented in the table indicate that this method seems to work for any $\theta$. We have added here the results for the case $n = 200$. It might surprise that the constant for "optimal" cv—bandwidths does not

TABLE 1

*Performance of MD and PL: Means (first line), standard deviations (second line) and mean squared error (third line) of $\widehat{\theta}$ for different $\theta_o$, models [see (14)], and bandwidths $h_\alpha = h_0 n^{-1/5}$, $\alpha = 1, 2$, for sample size $n = 100$. All numbers are calculated from 500 replications*

| | Both methods when using BF | | | | | | | | | | | | | | | | | |
|---|---|---|---|---|---|---|---|---|---|---|---|---|---|---|---|---|---|---|
| | **MD** | | | | | | | | | **PL** | | | | | | | | |
| $\theta_o$ | 0.00 | 0.50 | 1.0 | 0.00 | 0.50 | 1.0 | 0.00 | 0.50 | 1.0 | 0.00 | 0.50 | 1.0 | 0.00 | 0.50 | 1.0 | 0.00 | 0.50 | 1.0 |
| $h_0$ | | 0.3 | | | 0.4 | | | 0.5 | | | 0.2 | | | 0.3 | | | 0.4 | |
| Model 1 | | | | | | | | | | | | | | | | | | |
| | 0.02 | 0.53 | 0.92 | 0.02 | 0.53 | 0.92 | 0.03 | 0.56 | 0.92 | −0.00 | 0.43 | 0.83 | −0.01 | 0.43 | 0.83 | −0.00 | 0.43 | 0.83 |
| | 0.11 | 0.40 | 0.55 | 0.12 | 0.42 | 0.58 | 0.12 | 0.44 | 0.58 | 0.07 | 0.28 | 0.44 | 0.08 | 0.29 | 0.47 | 0.08 | 0.31 | 0.49 |
| | 0.01 | 0.16 | 0.31 | 0.01 | 0.18 | 0.34 | 0.02 | 0.20 | 0.34 | 0.01 | 0.08 | 0.22 | 0.01 | 0.09 | 0.24 | 0.01 | 0.10 | 0.27 |
| Model 2 | | | | | | | | | | | | | | | | | | |
| | 0.03 | 0.57 | 0.94 | 0.03 | 0.58 | 0.94 | 0.04 | 0.60 | 0.94 | −0.00 | 0.45 | 0.87 | −0.00 | 0.44 | 0.85 | −0.00 | 0.45 | 0.84 |
| | 0.15 | 0.44 | 0.56 | 0.16 | 0.46 | 0.57 | 0.16 | 0.47 | 0.58 | 0.01 | 0.31 | 0.46 | 0.10 | 0.32 | 0.47 | 0.10 | 0.33 | 0.50 |
| | 0.02 | 0.20 | 0.31 | 0.03 | 0.22 | 0.33 | 0.03 | 0.23 | 0.34 | 0.01 | 0.10 | 0.23 | 0.01 | 0.10 | 0.25 | 0.01 | 0.11 | 0.27 |
| Model 3 | | | | | | | | | | | | | | | | | | |
| | 0.05 | 0.60 | 0.96 | 0.07 | 0.61 | 0.96 | 0.08 | 0.63 | 0.97 | 0.00 | 0.46 | 0.87 | 0.00 | 0.45 | 0.86 | 0.00 | 0.45 | 0.86 |
| | 0.23 | 0.47 | 0.54 | 0.24 | 0.49 | 0.57 | 0.24 | 0.50 | 0.58 | 0.15 | 0.34 | 0.46 | 0.16 | 0.36 | 0.48 | 0.16 | 0.36 | 0.49 |
| | 0.05 | 0.23 | 0.29 | 0.06 | 0.25 | 0.33 | 0.07 | 0.27 | 0.34 | 0.02 | 0.12 | 0.23 | 0.02 | 0.13 | 0.26 | 0.02 | 0.13 | 0.26 |

Table 2

*Simulation results for different sample sizes n with cross validation bandwidth to minimize* (10) *with respect to θ. Numbers are calculated from 100 replications*

| **MD with cv-bandwidth** | | | | | |
|---|---|---|---|---|---|
| **$n$** | **100** | | | **200** | | |
| $\theta_o$ | **mean($\widehat{\theta}$)** | **std($\widehat{\theta}$)** | **mse** | **mean($\widehat{\theta}$)** | **std($\widehat{\theta}$)** | **mse** |
| 0.0 | 0.01 | 0.14 | 0.02 | 0.02 | 0.06 | 0.01 |
| 0.5 | 0.50 | 0.53 | 0.28 | 0.55 | 0.29 | 0.09 |
| 1.0 | 0.83 | 0.61 | 0.40 | 1.0 | 0.37 | 0.14 |

only change with $\theta$, but even more with $n$ (not shown in table). Have in mind that in small samples the second order terms of bias and variance are still quite influential and, thus, the rate $n^{-1/5}$ is to be taken carefully; compare with the above convergence-rate study.

A disadvantage of this cross validation procedure is that it is computationally rather expensive, and often rather hard to implement in practice. This is especially true if one wants to combine the cross validation method with the PL method. Sperlich, Linton and Van Keilegom (2007) discuss some alternative approaches like choosing $\theta$ and the bandwidth, simultaneously minimizing, respectively maximizing, the considered criteria function (8), respectively (10). In the same work are given results on the performance of the suggested bootstrap procedures which turn out there to perform reasonably well.

5.3. *Comparison with existing methods.* To our knowledge, the only existing method comparable to ours has been proposed by Linton, Chen, Wang and Härdle (1997). They considered the criterion functions

$$Q_3 = (\epsilon_\theta^T Z \ W \ Z^T \epsilon_\theta) \quad \text{and} \quad Q_4 = \frac{1}{n} \sum_{i=1}^{n} J_\theta(Y_i) - \ln\left\{\frac{1}{n}\epsilon_\theta^T \epsilon_\theta\right\},$$

where $\epsilon_\theta = (\epsilon_\theta^1, \ldots, \epsilon_\theta^n)^\top$ is the vector of residuals of the transformed model using $\theta$, while $Z = (Z_1, \ldots Z_n)^T$ are i.i.d. instruments with the property $E[Z_i \epsilon_\theta^i] = 0$. Here, $W$ is any symmetric positive definite weighting matrix, and $J_\theta$ is the Jacobian of the transformation $\Lambda_\theta$. When we tried to estimate $\theta$ in our simulation model (14), both criteria gave us always $-0.25$ for any data generating $\theta_o$. This was true for whichever smoother we used [in their article they just work with the marginal integration estimator]. The problem could come from the fact that they do not take care for the change of the total variation when transforming the response variable $Y$. Therefore, we have tried some modifications norming the criteria function by the total variation. Then the results change a lot, but still fail in estimating $\theta$.

## APPENDIX A: PROFILE LIKELIHOOD ESTIMATOR

To prove the asymptotic normality of the profile likelihood estimator, we will use Theorems 1 and 2 of Chen, Linton and Van Keilegom (2003) [abbreviated by CLV (2003) in the sequel]. Therefore, we need to define the space to which the nuisance function $s = (m, r, f, g, h)$ belongs. We define this space by $\mathcal{H}_{\mathrm{PL}} = \mathcal{M}^2 \times C_1^1(\mathbb{R})^3$, where $C_a^b(R)$ $(0 < a < \infty,\ 0 < b \leq 1,\ R \subset \mathbb{R}^k$ for some $k)$ is the set of all continuous functions $f : R \to \mathbb{R}$ for which

$$\sup_y |f(y)| + \sup_{y,y'} \frac{|f(y) - f(y')|}{|y - y'|^b} \leq a,$$

and where the space $\mathcal{M}$ depends on the model at hand. For instance, when the model is additive, a good choice for $\mathcal{M}$ is $\mathcal{M} = \sum_{\alpha=1}^d C_1^1(R_{X_\alpha})$, and when the model is multiplicative, $\mathcal{M} = \prod_{\alpha=1}^d C_1^1(R_{X_\alpha})$. We also need to define, according to CLV (2003), a norm for the space $\mathcal{H}_{\mathrm{PL}}$. Let

$$\|s\|_{\mathrm{PL}} = \sup_{\theta \in \Theta} \max\{\|m_\theta\|_\infty, \|r_\theta\|_\infty, \|f_\theta\|_2, \|g_\theta\|_2, \|h_\theta\|_2\},$$

where $\|\cdot\|_\infty$ ($\|\cdot\|_2$) denotes the $L_\infty$ ($L_2$) norm. Finally, let's denote $\|\cdot\|$ for the Euclidean norm.

We assume that the estimator $\widehat{m}_\theta$ is constructed based on a kernel function of degree $q_1$, which we assume of the form $K_1(u_1) \times \cdots \times K_1(u_d)$, and a bandwidth $h$. The required conditions on $K_1, q_1$ and $h$ are mentioned in the list of regularity conditions given below.

**A.1. Assumptions.** We assume throughout this appendix that the conditions stated below are satisfied. Condition A.1–A.7 are regularity conditions on the kernels, bandwidths, distributions $F_X$, $F_\varepsilon$, etc., whereas condition A.8 contains primitive conditions on the estimator $\widehat{m}_\theta$ that need to be checked depending on which model structure and which estimator $\widehat{m}_\theta$ one has chosen.

A.1 The probability density function $K_j$ $(j = 1, 2)$ is symmetric and has compact support, $\int u^k K_j(u)\, du = 0$ for $k = 1, \ldots, q_j - 1$, $\int u^{q_j} K_j(u)\, du \neq 0$ and $K_j$ is twice continuously differentiable.

A.2 $nh \to \infty$, $nh^{2q_1} \to 0$, $ng^6 (\log g^{-1})^{-2} \to \infty$ and $ng^{2q_2} \to 0$, where $q_1$ and $q_2$ are defined in condition A.1 and $q_1, q_2 \geq 4$.

A.3 The density $f_X$ is bounded away from zero and infinity and is Lipschitz continuous on the compact support $\mathcal{X}$.

A.4 The functions $m_\theta(x)$ and $\dot{m}_\theta(x)$ are $q_1$ times continuously differentiable with respect to the components of $x$ on $\mathcal{X} \times \mathcal{N}(\theta_o)$, and all derivatives up to order $q_1$ are bounded, uniformly in $(x, \theta)$ in $\mathcal{X} \times \mathcal{N}(\theta_o)$.

A.5 The transformation $\Lambda_\theta(y)$ is three times continuously differentiable in both $\theta$ and $y$, and there exists a $\delta > 0$ such that

$$E\left[ \sup_{\|\theta' - \theta\| \leq \delta} \left| \frac{\partial^{k+l}}{\partial y^k\, \partial \theta^l} \Lambda_{\theta'}(Y) \right| \right] < \infty$$

for all $\theta$ in $\Theta$ and all $0 \leq k + l \leq 3$.

A.6 The distribution $F_{\varepsilon(\theta)}(y)$ is three times continuously differentiable with respect to $y$ and $\theta$, and

$$\sup_{\theta, y} \left| \frac{\partial^{k+l}}{\partial y^k \, \partial \theta^l} F_{\varepsilon(\theta)}(y) \right| < \infty$$

for all $0 \leq k + l \leq 2$.

A.7 For all $\eta > 0$, there exists $\epsilon(\eta) > 0$ such that

$$\inf_{\|\theta - \theta_o\| > \eta} \|G_{\mathrm{PL}}(\theta, s_\theta)\| \geq \epsilon(\eta) > 0.$$

Moreover, the matrix $\Gamma_{1\mathrm{PL}}$ is of full (column) rank.

A.8 The estimators $\widehat{m}_o$ and $\dot{\widehat{m}}_o$ can be written as

$$\widehat{m}_o(x) - m_o(x) = \frac{1}{nh} \sum_{i=1}^{n} \sum_{\alpha=1}^{d} K_1 \left( \frac{x_\alpha - X_{\alpha i}}{h} \right) v_{o1\alpha}(X_{\alpha i}, \varepsilon_i)$$
$$+ \frac{1}{n} \sum_{i=1}^{n} v_{o2}(X_i, \varepsilon_i) + \widehat{v}_o(x)$$

and

$$\dot{\widehat{m}}_o(x) - \dot{m}_o(x) = \frac{1}{nh} \sum_{i=1}^{n} \sum_{\alpha=1}^{d} K_1 \left( \frac{x_\alpha - X_{\alpha i}}{h} \right) w_{o1\alpha}(X_{\alpha i}, \varepsilon_i)$$
$$+ \frac{1}{n} \sum_{i=1}^{n} w_{o2}(X_i, \varepsilon_i) + \widehat{w}_o(x),$$

where $\sup_x |\widehat{v}_o(x)| = o_p(n^{-1/2})$, $\sup_x |\widehat{w}_o(x)| = o_p(n^{-1/2})$, the functions $v_{o1\alpha}(x, e)$ and $w_{o1\alpha}(x, e)$ are $q_1$ times continuously differentiable with respect to the components of $x$, their derivatives up to order $q_1$ are bounded, uniformly in $x$ and $e$, $E(v_{o2}(X, \varepsilon)) = 0$ and $E(w_{o2}(X, \varepsilon)) = 0$. Moreover, with probability tending to 1, $\widehat{m}_\theta, \dot{\widehat{m}}_\theta \in \mathcal{M}$, $\sup_{\theta \in \Theta} \|\widehat{m}_\theta - m_\theta\| = o_p(1)$, $\sup_{\theta \in \Theta} \|\dot{\widehat{m}}_\theta - \dot{m}_\theta\| = o_p(1)$, $\|\widehat{m}_\theta - m_\theta\| = o_p(n^{-1/4})$ and $\|\dot{\widehat{m}}_\theta - \dot{m}_\theta\| = o_p(n^{-1/4})$ uniformly over all $\theta$ with $\|\theta - \theta_o\| = o(1)$, and

$$\sup_x |(\widehat{m}_\theta - \dot{m}_\theta)(x) - (\dot{\widehat{m}}_\theta - \dot{m}_o)(x)| = o_p(1)\|\theta - \theta_o\| + O_p(n^{-1/2})$$

for all $\theta$ with $\|\theta - \theta_o\| = o(1)$. Finally, the space $\mathcal{M}$ satisfies $\int \sqrt{\log N(\lambda, \mathcal{M}, \|\cdot\|_\infty)} \, d\lambda < \infty$, where $N(\lambda, \mathcal{M}, \|\cdot\|_\infty)$ is the covering number with respect to the norm $\|\cdot\|_\infty$ of the class $\mathcal{M}$, that is, the minimal number of balls of $\|\cdot\|_\infty$-radius $\lambda$ needed to cover $\mathcal{M}$.

**A.2. Proof of Theorem 4.1.** The proof consists of verifying the conditions given in Theorem 1 (regarding consistency) and Theorem 2 (regarding asymptotic normality) in CLV (2003). In Lemmas A.4–A.11 below, we verify these conditions. The result then follows immediately from those lemmas, assuming that the primitive conditions on $\widehat{m}_\theta$ and the regularity conditions stated in A.1–A.8 hold true. Before checking the conditions of these theorems, we first need to show three preliminary Lemmas A.1–A.3 which give asymptotic expansions for the estimators $f_\varepsilon$, $\hat{f}'_\varepsilon$ and $\hat{\dot{f}}_\varepsilon$. The proofs of all lemmas are deferred to Section A.3.

Lemma A.1. *For all $y \in \mathbb{R}$,*

$$\widehat{f}_\varepsilon(y) - f_\varepsilon(y) = n^{-1} \sum_{i=1}^n K_{2g}(\varepsilon_i - y) - f_\varepsilon(y)$$

$$+ f'_\varepsilon(y)\, n^{-1} \sum_{i=1}^n \left[ \sum_{\alpha=1}^d v_{o1\alpha}(X_{\alpha i}, \varepsilon_i) f_{X_\alpha}(X_{\alpha i}) + v_{o2}(\varepsilon_i) \right]$$

$$+ \widehat{r}_o(y),$$

*where $\sup_y |\widehat{r}_o(y)| = o_p(n^{-1/2})$, and where the functions $v_{o1\alpha}$ and $v_{o2}$ are defined in Assumption A.8. Moreover,*

$$\sup_y \sup_{\theta \in \Theta} |\widehat{f}_{\varepsilon(\theta)}(y) - f_{\varepsilon(\theta)}(y)| = o_p(1)$$

*and*

$$\sup_y \sup_{\|\theta - \theta_o\| \le \delta_n} |\widehat{f}_{\varepsilon(\theta)}(y) - f_{\varepsilon(\theta)}(y)| = o_p(n^{-1/4})$$

*for all $\delta_n = o(1)$.*

In a similar way as for Lemma A.1, we can prove the following two results. The proofs are omitted.

Lemma A.2. *For all $y \in \mathbb{R}$,*

$$\hat{\dot{f}}_\varepsilon(y) - \dot{f}_\varepsilon(y) = (ng)^{-1} \sum_{i=1}^n K'_{2g}(\varepsilon_i - y)(\dot{\Lambda}_\theta(Y_i) - \dot{m}_\theta(X_i)) - \dot{f}_\varepsilon(y)$$

$$+ \dot{f}'_\varepsilon(y)\, n^{-1} \sum_{i=1}^n \left[ \sum_{\alpha=1}^d v_{o1\alpha}(X_{\alpha i}, \varepsilon_i) f_{X_\alpha}(X_{\alpha i}) + v_{o2}(\varepsilon_i) \right]$$

$$+ f'_\varepsilon(y)\, n^{-1} \sum_{i=1}^n \left[ \sum_{\alpha=1}^d w_{o1\alpha}(X_{\alpha i}, \varepsilon_i) f_{X_\alpha}(X_{\alpha i}) + w_{o2}(\varepsilon_i) \right]$$

$$+ \widehat{r}_o(y),$$

*where* $\sup_{\theta,y} |\widehat{r}_o(y)| = o_p(n^{-1/2})$. *Moreover,*

$$\sup_y \sup_{\theta \in \Theta} |\dot{\widehat{f}}_{\varepsilon(\theta)}(y) - \dot{f}_{\varepsilon(\theta)}(y)| = o_p(1)$$

*and*

$$\sup_y \sup_{\|\theta - \theta_o\| \le \delta_n} |\dot{\widehat{f}}_{\varepsilon(\theta)}(y) - \dot{f}_{\varepsilon(\theta)}(y)| = o_p(n^{-1/4})$$

*for all* $\delta_n = o(1)$.

LEMMA A.3.   *For all* $y \in \mathbb{R}$,

$$\widehat{f}'_{\varepsilon}(y) - f'_{\varepsilon}(y) = (ng)^{-1} \sum_{i=1}^{n} K'_{2g}(y - \varepsilon_i) - f'_{\varepsilon}(y)$$

$$+ f''_{\varepsilon}(y) \, n^{-1} \sum_{i=1}^{n} \left[ \sum_{\alpha=1}^{d} v_{o1\alpha}(X_{\alpha i}, \varepsilon_i) f_{X_\alpha}(X_{\alpha i}) + v_{o2}(\varepsilon_i) \right]$$

$$+ \widehat{r}_o(y),$$

*where* $\sup_y |\widehat{r}_o(y)| = o_p(n^{-1/2})$. *Moreover,*

$$\sup_y \sup_{\theta \in \Theta} |\widehat{f}'_{\varepsilon(\theta)}(y) - f'_{\varepsilon(\theta)}(y)| = o_p(1)$$

*and*

$$\sup_y \sup_{\|\theta - \theta_o\| \le \delta_n} |\widehat{f}'_{\varepsilon(\theta)}(y) - f'_{\varepsilon(\theta)}(y)| = o_p(n^{-1/4})$$

*for all* $\delta_n = o(1)$.

LEMMA A.4.   *Uniformly for all* $\theta \in \Theta$, $G_{\mathrm{PL}}(\theta, s)$ *is continuous (with respect to the* $\| \cdot \|_{\mathrm{PL}}$-*norm) in* $s$ *at* $s = s_\theta$.

LEMMA A.5.

$$\sup_y \sup_{\theta \in \Theta} |\widehat{f}_{\varepsilon(\theta)}(y) - f_{\varepsilon(\theta)}(y)| = o_p(1),$$

$$\sup_y \sup_{\theta \in \Theta} |\dot{\widehat{f}}_{\varepsilon(\theta)}(y) - \dot{f}_{\varepsilon(\theta)}(y)| = o_p(1)$$

*and*

$$\sup_y \sup_{\theta \in \Theta} |\widehat{f}'_{\varepsilon(\theta)}(y) - f'_{\varepsilon(\theta)}(y)| = o_p(1).$$

Lemma A.6.  *For all sequences of positive numbers* $\delta_n = o(1)$,

$$\sup_{\theta \in \Theta, \|s - s_\theta\|_{\mathrm{PL}} \leq \delta_n} \|G_{n\mathrm{PL}}(\theta, s) - G_{\mathrm{PL}}(\theta, s)\| = o_p(1).$$

Lemma A.7.  *The ordinary partial derivative in* $\theta$ *of* $G_{\mathrm{PL}}(\theta, s_\theta)$, *denoted* $\Gamma_{1\mathrm{PL}}(\theta, s_\theta)$, *exists in a neighborhood of* $\theta_o$, *is continuous at* $\theta = \theta_o$, *and the matrix* $\Gamma_{1\mathrm{PL}} = \Gamma_{1\mathrm{PL}}(\theta_o, s_o)$ *is of full (column) rank.*

For any $\theta \in \Theta$, we say that $G_{\mathrm{PL}}(\theta, s)$ is pathwise differentiable at $s$ in the direction $[\overline{s} - s]$ if $\{s + \tau(\overline{s} - s) : \tau \in [0, 1]\} \subset \mathcal{H}_{\mathrm{PL}}$ and $\lim_{\tau \to 0}[G_{\mathrm{PL}}(\theta, s + \tau(\overline{s} - s)) - G_{\mathrm{PL}}(\theta, s)]/\tau$ exists; we denote the limit by $\Gamma_{2\mathrm{PL}}(\theta, s)[\overline{s} - s]$.

Lemma A.8.  *The pathwise derivative* $\Gamma_{2\mathrm{PL}}(\theta, s_\theta)$ *of* $G_{\mathrm{PL}}(\theta, s_\theta)$ *exists in all directions* $s - s_\theta$ *and satisfies the following:*

(i)  $\|G_{\mathrm{PL}}(\theta, s) - G_{\mathrm{PL}}(\theta, s_\theta) - \Gamma_{2\mathrm{PL}}(\theta, s_\theta)[s - s_\theta]\| \leq c\|s - s_\theta\|_{\mathrm{PL}}^2$

*for all* $\theta$ *with* $\|\theta - \theta_o\| = o(1)$, *all* $s$ *with* $\|s - s_\theta\|_{\mathrm{PL}} = o(1)$, *some constant* $c < \infty$;

(ii)  $\|\Gamma_{2\mathrm{PL}}(\theta, s_\theta)[\widehat{s}_\theta - s_\theta] - \Gamma_{2\mathrm{PL}}(\theta_o, s_o)[\widehat{s}_o - s_o]\|$

$$\leq c\|\theta - \theta_o\| \times o_p(1) + O_p(n^{-1/2})$$

*for all* $\theta$ *with* $\|\theta - \theta_o\| = o(1)$, *where* $\widehat{s} = (\widehat{m}, \dot{\widehat{m}}, \widehat{f}_\varepsilon, \dot{\widehat{f}}_\varepsilon, \widehat{f}'_\varepsilon)$.

Lemma A.9.  *With probability tending to one,* $\widehat{f}_\varepsilon, \dot{\widehat{f}}_\varepsilon, \widehat{f}'_\varepsilon \in C_1^1(\mathbb{R})$. *More-over,*

$$\sup_y \sup_{\|\theta - \theta_o\| \leq \delta_n} |\widehat{f}_{\varepsilon(\theta)}(y) - f_{\varepsilon(\theta)}(y)| = o_p(n^{-1/4}),$$

$$\sup_y \sup_{\|\theta - \theta_o\| \leq \delta_n} |\dot{\widehat{f}}_{\varepsilon(\theta)}(y) - \dot{f}_{\varepsilon(\theta)}(y)| = o_p(n^{-1/4})$$

*and*

$$\sup_y \sup_{\|\theta - \theta_o\| \leq \delta_n} |\widehat{f}'_{\varepsilon(\theta)}(y) - f'_{\varepsilon(\theta)}(y)| = o_p(n^{-1/4}),$$

*for any* $\delta_n = o(1)$.

Lemma A.10.  *For all sequences of positive numbers* $\{\delta_n\}$ *with* $\delta_n = o(1)$,

$$\sup_{\|\theta - \theta_o\| \leq \delta_n, \|s - s_\theta\|_{\mathrm{PL}} \leq \delta_n} \|G_{n\mathrm{PL}}(\theta, s) - G_{\mathrm{PL}}(\theta, s) - G_{n\mathrm{PL}}(\theta_o, s_o)\| = o_p(n^{-1/2}).$$

Lemma A.11.

$$\sqrt{n}\{G_{n\mathrm{PL}}(\theta_o, s_o) + \Gamma_{2\mathrm{PL}}(\theta_o, s_o)[\widehat{s} - s_o]\} \Longrightarrow N(0, \mathrm{Var}\{G_{1\mathrm{PL}}(\theta_o, s_o)\}).$$

**A.3. Proofs of Lemmas A.1–A.11.**

PROOF OF LEMMA A.1.    Write

$$\widehat{f}_\varepsilon(y) - f_\varepsilon(y)$$

$$= \frac{1}{ng}\sum_{i=1}^n K'_{2g}(\varepsilon_i - y)(\widehat{\varepsilon}_i - \varepsilon_i)$$

$$+ \frac{1}{n}\sum_{i=1}^n K_{2g}(\varepsilon_i - y) - f_\varepsilon(y) + o_p(n^{-1/2})$$

(15)
$$= -\frac{1}{ng}\sum_{i=1}^n K'_{2g}(\varepsilon_i - y)\left\{\frac{1}{n}\sum_{k=1}^n\sum_{\alpha=1}^d K_{1h}(X_{\alpha i} - X_{\alpha k})v_{o1\alpha}(X_{\alpha k}, \varepsilon_k)\right.$$

$$\left. + \frac{1}{n}\sum_{k=1}^n v_{o2}(\varepsilon_k) + \widehat{v}_o(X_i)\right\}$$

$$+ \frac{1}{n}\sum_{i=1}^n K_{2g}(\varepsilon_i - y) - f_\varepsilon(y) + o_p(n^{-1/2})$$

$$= \frac{1}{n^2}\sum_{\alpha=1}^d\sum_{i,k=1}^n v_{o1\alpha}(X_{\alpha k}, \varepsilon_k)\varphi_{nik} + f'_\varepsilon(y)\frac{1}{n}\sum_{k=1}^n v_{o2}(\varepsilon_k)$$

(16)
$$+ \frac{1}{n}\sum_{i=1}^n K_{2g}(\varepsilon_i - y) - f_\varepsilon(y) + o_p(n^{-1/2}),$$

where $\varphi_{nik} = -\frac{1}{g}K'_{2g}(\varepsilon_i - y)K_{1h}(X_{\alpha i} - X_{\alpha k})$. Since $E(\varphi_{nik}|X_k) = f'_\varepsilon(y)f_{X_\alpha}(X_{\alpha k}) + o_p(1)$, it follows that (16) equals

$$f'_\varepsilon(y)\frac{1}{n}\sum_{k=1}^n\left[\sum_{\alpha=1}^d v_{o1\alpha}(X_{\alpha k}, \varepsilon_k)f_{X_\alpha}(X_{\alpha k}) + v_{o2}(\varepsilon_k)\right]$$

$$+ \frac{1}{n}\sum_{i=1}^n K_{2g}(\varepsilon_i - y) - f_\varepsilon(y) + o_p(n^{-1/2}). \qquad \square$$

PROOF OF LEMMA A.4.    Note that

$$G_{\mathrm{PL}}(\theta, s)$$

$$= E\left[\frac{1}{f(\varepsilon(\theta, m))}\{g(\varepsilon(\theta, m))(\dot{\Lambda}_\theta(Y) - r(X)) + h(\varepsilon(\theta, m))\} + \frac{\dot{\Lambda}'_\theta(Y)}{\Lambda'_\theta(Y)}\right],$$

which is continuous in $s$ at $s = s_\theta$, provided conditions A.4–A.6 are satisfied.
$\square$

PROOF OF LEMMA A.5.  This follows from Lemmas A.1–A.3.  □

PROOF OF LEMMA A.6.  The proof is similar to (but easier than) that of Lemma A.10. We therefore omit the proof.  □

PROOF OF LEMMA A.7.  This follows from Assumption A.7.  □

PROOF OF LEMMA A.8.  Some straightforward calculations show that

$$
\begin{aligned}
\Gamma_{2\mathrm{PL}}&(\theta, s_\theta)[\widehat{s}_\theta - s_\theta] \\
&= \lim_{\tau \to 0} \frac{1}{\tau}\{G_{\mathrm{PL}}(\theta, s_\theta + \tau(\widehat{s}_\theta - s_\theta)) - G_{\mathrm{PL}}(\theta, s_\theta)\} \\
&= E\bigg[\bigg\{\frac{f'_{\varepsilon(\theta)}(\varepsilon_\theta)}{f^2_{\varepsilon(\theta)}(\varepsilon_\theta)}(\widehat{m}_\theta - m_\theta)(X) - \frac{(\widehat{f}_{\varepsilon(\theta)} - f_{\varepsilon(\theta)})(\varepsilon_\theta)}{f^2_{\varepsilon(\theta)}(\varepsilon_\theta)}\bigg\} \\
&\qquad \times \bigg\{f'_{\varepsilon(\theta)}(\varepsilon_\theta)[\dot{\Lambda}_\theta(Y) - \dot{m}_\theta(X)] + \dot{f}_{\varepsilon(\theta)}(\varepsilon_\theta)\bigg\} \\
&\quad + \frac{1}{f_{\varepsilon(\theta)}(\varepsilon_\theta)}\{-f''_{\varepsilon(\theta)}(\varepsilon_\theta)[\dot{\Lambda}_\theta(Y) - \dot{m}_\theta(X)](\widehat{m}_\theta - m_\theta)(X) \\
&\qquad\qquad + (\widehat{f}_{\varepsilon(\theta)} - f'_{\varepsilon(\theta)})(\varepsilon_\theta)[\dot{\Lambda}_\theta(Y) - \dot{m}_\theta(X)] \\
&\qquad\qquad - f'_{\varepsilon(\theta)}(\varepsilon_\theta)(\widehat{\dot{m}}_\theta - \dot{m}_\theta)(X) \\
&\qquad\qquad + (\widehat{\dot{f}}_{\varepsilon(\theta)} - \dot{f}_{\varepsilon(\theta)})(\varepsilon_\theta) - \dot{f}'_{\varepsilon(\theta)}(\varepsilon_\theta)(\widehat{m}_\theta - m_\theta)(X)\}\bigg].
\end{aligned}
\tag{17}
$$

The first part of Lemma A.8 now follows immediately. The second part follows from the uniform consistency of $\widehat{m}$, $\widehat{\dot{m}}$, $\widehat{f}_{\varepsilon(\theta)}$, $\widehat{\dot{f}}_{\varepsilon(\theta)}$ and $\widehat{f}'_{\varepsilon(\theta)}$, and from the fact that

$$
\sup_x |(\widehat{\dot{m}}_\theta - \dot{m}_\theta)(x) - (\widehat{\dot{m}}_o - \dot{m}_o)(x)| = o_p(1)\|\theta - \theta_o\| + O_p(n^{-1/2}),
$$

which follows from Assumption A.8.  □

PROOF OF LEMMA A.9.  This follows from Lemmas A.1–A.3.  □

PROOF OF LEMMA A.10.  We will make use of Theorem 3 in Chen, Linton and Van Keilegom (2003). According to this result, we need to prove that

(i)

$$
E\bigg[\sup_{\|\theta' - \theta\| < \eta, \|s' - s\|_{\mathrm{PL}} < \eta} |g_{\mathrm{PL}}(X, Y, \theta', s') - g_{\mathrm{PL}}(X, Y, \theta, s)|^2\bigg] \le K\eta^2,
$$

for all $(\theta, s) \in \Theta \times \mathcal{H}_{\mathrm{PL}}$, all $\eta > 0$ and for some $K > 0$.

(ii)

$$\int_0^\infty \sqrt{\log N(\lambda, \mathcal{H}_{\mathrm{PL}}, \|\cdot\|_{\mathrm{PL}})}\, d\lambda < \infty.$$

Part (ii) follows from Corollary 2.7.4 in van der Vaart and Wellner (1996), together with Assumption A.8. Part (i) follows from the mean value theorem, together with the differentiability conditions imposed on the functions of which the function $g_{\mathrm{PL}}$ is composed. $\square$

PROOF OF LEMMA A.11. Combining the formula of $\Gamma_{2\mathrm{PL}}(\theta_o, s_o)$ given in (17) with the representations of $\widehat{f}_{\varepsilon(\theta)}$, $\widehat{\dot{f}}_{\varepsilon(\theta)}$ and $\widehat{f}'_{\varepsilon(\theta)}$ given in Lemmas A.1–A.3, we obtain after some calculations

$$G_{n\mathrm{PL}}(\theta_o, s_o) + \Gamma_{2\mathrm{PL}}(\theta_o, s_o)[\widehat{s} - s_o]$$

$$= n^{-1} \sum_{i=1}^n \left\{ \frac{1}{f_\varepsilon(\varepsilon_i)} [f'_\varepsilon(\varepsilon_i)\{\dot{\Lambda}_o(Y_i) - \dot{m}_o(X_i)\} + \dot{f}_\varepsilon(\varepsilon_i)] + \frac{\dot{\Lambda}'_o(Y_i)}{\Lambda'_o(Y_i)} \right\}$$

$$+ E\left[ -\frac{1}{f_\varepsilon^2(\varepsilon)} \left\{ \frac{1}{ng} \sum_{i=1}^n K_2\left(\frac{\varepsilon_i - \varepsilon}{g}\right) - f_\varepsilon(\varepsilon) \right\} \right.$$

(18)

$$\times \{f'_\varepsilon(\varepsilon)[\dot{\Lambda}_o(Y) - \dot{m}_o(X)] + \dot{f}_\varepsilon(\varepsilon)\}$$

$$+ \frac{1}{f_\varepsilon(\varepsilon)} \left\{ -\frac{1}{ng^2} \sum_{i=1}^n K'_2\left(\frac{\varepsilon_i - \varepsilon}{g}\right) - f'_\varepsilon(\varepsilon) \right\} \{\dot{\Lambda}_o(Y) - \dot{m}_o(X)\}$$

$$+ \frac{1}{f_\varepsilon(\varepsilon)} \left\{ \frac{1}{ng^2} \sum_{i=1}^n K'_2\left(\frac{\varepsilon_i - \varepsilon}{g}\right)(\dot{\Lambda}_o(Y_i) - \dot{m}_o(X_i)) - \dot{f}_\varepsilon(\varepsilon) \right\} \right]$$

$$+ o_p(n^{-1/2}).$$

We next show that

(19)
$$E\left[ \frac{\dot{f}_\varepsilon(\varepsilon)}{f_\varepsilon(\varepsilon)} \right] = 0,$$

(20) $\quad E\left[ \dfrac{1}{f_\varepsilon(\varepsilon)} \left\{ \dfrac{1}{ng^2} \displaystyle\sum_{i=1}^n K'_2\left(\dfrac{\varepsilon_i - \varepsilon}{g}\right)(\dot{\Lambda}_o(Y_i) - \dot{m}_o(X_i)) - \dot{f}_\varepsilon(\varepsilon) \right\} \right] = 0$

and

$$E\left[ -\frac{1}{f_\varepsilon^2(\varepsilon)} \left\{ \frac{1}{ng} \sum_{i=1}^n K_2\left(\frac{\varepsilon_i - \varepsilon}{g}\right) \right\} \{f'_\varepsilon(\varepsilon)[\dot{\Lambda}_o(Y) - \dot{m}_o(X)] + \dot{f}_\varepsilon(\varepsilon)\} \right.$$

(21)

$$\left. + \frac{1}{f_\varepsilon(\varepsilon)} \left\{ -\frac{1}{ng^2} \sum_{i=1}^n K'_2\left(\frac{\varepsilon_i - \varepsilon}{g}\right) \right\} \{\dot{\Lambda}_o(Y) - \dot{m}_o(X)\} \right] = 0.$$

It then follows that only the first term on the right-hand side of (18) [i.e., the term $G_{n\mathrm{PL}}(\theta_o, s_o)$] is nonzero, from which the result follows. We start by showing (19):

$$E\left[\frac{\dot{f}_\varepsilon(\varepsilon)}{f_\varepsilon(\varepsilon)}\right] = \int \dot{f}_\varepsilon(y)\,dy = \frac{\partial}{\partial\theta}\int f_{\varepsilon(\theta)}(y)\,dy\bigg\downarrow_{\theta=\theta_o} = 0,$$

since $\int f_{\varepsilon(\theta)}(y)\,dy = 1$. Next, consider (20). The left-hand side equals

$$\frac{1}{ng^2}\sum_{i=1}^n (\dot{\Lambda}_o(Y_i) - \dot{m}_o(X_i))E\left[\frac{1}{f_\varepsilon(\varepsilon)}K_2'\left(\frac{\varepsilon_i - \varepsilon}{g}\right)\right] - E\left[\frac{\dot{f}_\varepsilon(\varepsilon)}{f_\varepsilon(\varepsilon)}\right]$$

$$= \frac{1}{ng}\sum_{i=1}^n (\dot{\Lambda}_o(Y_i) - \dot{m}_o(X_i))\int K_2'(u)\,du = 0.$$

Finally, for (22), note that the left-hand side can be written as

$$\frac{1}{ng}\sum_{i=1}^n E\left[\frac{1}{f_\varepsilon^2(\varepsilon)}\left\{-K_2\left(\frac{\varepsilon_i - \varepsilon}{g}\right)\frac{d}{d\theta}f_{\varepsilon(\theta)}(\varepsilon(\theta))\downarrow_{\theta=\theta_o}\right.\right.$$

$$\left.\left. +\frac{d}{d\theta}K_2\left(\frac{\varepsilon_i - \varepsilon(\theta)}{g}\right)\bigg\downarrow_{\theta=\theta_o}f_\varepsilon(\varepsilon)\right\}\right]$$

$$= \frac{1}{ng}\sum_{i=1}^n E\left[\frac{d}{d\theta}\frac{K_2((\varepsilon_i - \varepsilon(\theta))/g)}{f_{\varepsilon(\theta)}(\varepsilon(\theta))}\bigg\downarrow_{\theta=\theta_o}\right]$$

$$= \frac{1}{ng}\sum_{i=1}^n \frac{d}{d\theta}\int K_2\left(\frac{\varepsilon_i - e}{g}\right)de = 0,$$

since $\int K_2(\frac{\varepsilon_i - e}{g})\,de = g$. This finishes the proof. $\quad\square$

## APPENDIX B: MD ESTIMATOR

**B.1. Assumptions.** We assume throughout this appendix that Assumptions B.1–B.8 given below are valid.

B.1 The probability density function $K_1$ is symmetric and has compact support, $\int u^k K_1(u)\,du = 0$ for $k = 1, \ldots, q_1 - 1$, $\int u^{q_1} K_1(u)\,du \neq 0$ and $K_1$ is twice continuously differentiable.

B.2 $nh \to \infty$ and $nh^{2q_1} \to 0$, where $q_1$ is defined in condition B.1 and $q_1 \geq 4$.

B.3 The density $f_X$ is bounded away from zero and infinity and is Lipschitz continuous on the compact support $\mathcal{X}$.

B.4 The function $m_\theta(x)$ is $q_1$ times continuously differentiable with respect to the components of $x$ on $\mathcal{X} \times \mathcal{N}(\theta_o)$, and all derivatives up to order $q_1$ are bounded, uniformly in $(x, \theta)$ in $\mathcal{X} \times \mathcal{N}(\theta_o)$.

B.5 The transformation $\Lambda_\theta(y)$ is twice continuously differentiable in both $\theta$ and $y$, and there exists a $\delta > 0$ such that

$$E\left[\sup_{\|\theta - \theta'\| \le \delta} |\lambda_{\theta'}(Y)|^k\right] < \infty$$

for all $k$ and for all $\theta$ in $\Theta$.

B.6 The distribution $F_\varepsilon(y)$ is twice continuously differentiable with respect to $y$, and $\sup_y |f'_\varepsilon(y)| < \infty$.

B.7 For all $\eta > 0$, there exists $\epsilon(\eta) > 0$ such that

$$\inf_{\|\theta - \theta_o\| > \eta} \|G_{\mathrm{MD}}(\theta, m_\theta)\|_2 \ge \epsilon(\eta) > 0.$$

Moreover, the matrix $\Gamma_{1\mathrm{MD}}(x, e)$ (defined in Section 4) is of full (column) rank for a set of positive $\mu$-measure $(x, e)$.

B.8 The estimator $\widehat{m}_o$ can be written as

$$\widehat{m}_o(x) - m_o(x) = \frac{1}{nh} \sum_{i=1}^n \sum_{\alpha=1}^d K_1\left(\frac{x_\alpha - X_{\alpha i}}{h}\right) v_{o1\alpha}(X_{\alpha i}, \varepsilon_i)$$
$$+ \frac{1}{n} \sum_{i=1}^n v_{o2}(X_i, \varepsilon_i) + \widehat{v}_o(x),$$

where $\sup_x |\widehat{v}_o(x)| = o_p(n^{-1/2})$, the function $v_{o1\alpha}(x, e)$ is $q_1$ times continuously differentiable with respect to the components of $x$, their derivatives up to order $q_1$ are bounded, uniformly in $x$ and $e$, $E(v_{o2}(X, \varepsilon)) = 0$. Moreover, with probability tending to 1, $\widehat{m}_\theta \in \mathcal{M}$, $\sup_{\theta \in \Theta} \|\widehat{m}_\theta - m_\theta\| = o_p(1)$, $\|\widehat{m}_\theta - m_\theta\| = o_p(n^{-1/4})$ uniformly over all $\theta$ with $\|\theta - \theta_o\| = o(1)$, and

$$\sup_x |(\widehat{m}_\theta - m_\theta)(x) - (\widehat{m}_o - m_o)(x)| = o_p(1)\|\theta - \theta_o\| + O_p(n^{-1/2})$$

for all $\theta$ with $\|\theta - \theta_o\| = o(1)$. Finally, the space $\mathcal{M}$ satisfies $\int \sqrt{\log N(\lambda, \mathcal{M}, \|\cdot\|_\infty)}\, d\lambda < \infty$.

**B.2. Proof of Theorem 4.2.** We use a generalization of Theorems 1 (about consistency) and 2 (about asymptotic normality) of Chen, Linton and Van Keilegom (2003), henceforth, CLV (2003). Below, we state the primitive conditions under which these results are valid (see Lemmas B.1–B.6). Their proof is given in Section B.3.

Given these lemmas, we have the desired result. We just reprieve the last part of the argument because it is slightly different from CLV (2003) due to the different norm. Note that

$$F_{\varepsilon(\theta, m)}(e) = \Pr[\Lambda_\theta(Y) - m(X) \le e]$$

$$= \Pr[Y \leq \Lambda_\theta^{-1}(m(X) + e)]$$
$$= \Pr[\varepsilon \leq \Lambda_o(\Lambda_\theta^{-1}(m(X) + e)) - m_o(X)]$$
$$= E F_\varepsilon[\Lambda_o(\Lambda_\theta^{-1}(m(X) + e)) - m_o(X)].$$

Likewise, $F_{X,\varepsilon(\theta,m)}$ satisfies

$$F_{X,\varepsilon(\theta,m)}(x,e) = \Pr[X \leq x, \Lambda_\theta(Y) - m(X) \leq e]$$
$$= E \Pr[X \leq x, \varepsilon \leq \Lambda_o(\Lambda_\theta^{-1}(m(X) + e)) - m_o(X)]$$
$$= E[1(X \leq x) F_\varepsilon[\Lambda_o(\Lambda_\theta^{-1}(m(X) + e)) - m_o(X)]].$$

Define

$$G_{\mathrm{MD}}(\theta,m)(x,e) = F_{X,\varepsilon(\theta,m)}(x,e) - F_X(x) F_{\varepsilon(\theta,m)}(e).$$

Define now the stochastic processes

$$L_n(x,e) = \sqrt{n}[\widehat{F}_{X,\varepsilon}(x,e) - F_{X,\varepsilon}(x,e)]$$
$$- F_X(x)\sqrt{n}[\widehat{F}_\varepsilon(e) - F_\varepsilon(e)] - F_\varepsilon(e)\sqrt{n}[\widehat{F}_X(x) - F_X(x)]$$

and

$$\mathcal{L}_n(\theta)(x,e) = L_n(x,e) + \Gamma_{\mathrm{1MD}}(x,e)(\theta - \theta_o) + [\Gamma_{\mathrm{2MD}}(\theta_o,m_o)(\widehat{m} - m_o)](x,e),$$

where for any $\theta \in \Theta$ and any $m,\overline{m} \in \mathcal{M}$, $\Gamma_{\mathrm{2MD}}(\theta,m)(\overline{m} - m)(x,e)$ is defined in the following way. We say that $G_{\mathrm{MD}}(\theta,m)$ is pathwise differentiable at $m$ in the direction $[\overline{m} - m]$ at $(x,e)$ if $\{m + \tau(\overline{m} - m) : \tau \in [0,1]\} \subset \mathcal{M}$ and $\lim_{\tau \to 0}[G_{\mathrm{MD}}(\theta, m + \tau(\overline{m} - m))(x,e) - G_{\mathrm{MD}}(\theta,m)(x,e)]/\tau$ exist; we denote the limit by $\Gamma_{\mathrm{2MD}}(\theta,m)[\overline{m} - m](x,e)$.

A consequence of Lemmas B.1–B.6 is that

$$\sup_{\|\theta - \theta_o\| \leq \delta_n} \|G_{n\mathrm{MD}}(\theta, \widehat{m}_\theta) - \mathcal{L}_n(\theta)\|_2^2 = o_p(n^{-1/2}),$$

which means we can effectively deal with the minimizer of $\mathcal{L}_n(\theta)$, say, $\overline{\theta}$. Note that $\overline{\theta}$ has an explicit solution and, indeed,

$$\sqrt{n}(\overline{\theta} - \theta_o) = -\left[\int \Gamma_{\mathrm{1MD}} \Gamma_{\mathrm{1MD}}^\top(x,e)\, d\mu(x,e)\right]^{-1}$$
$$\times \int [L_n(x,e) + [\Gamma_{\mathrm{2MD}}(\theta_o,m_o)(\widehat{m} - m_o)](x,e)]$$
$$\times \Gamma_{\mathrm{1MD}}(x,e)\, d\mu(x,e).$$

Then apply Lemma B.6 below to get the desired result.

LEMMA B.1. *Uniformly for all $\theta \in \Theta$, $G_{\mathrm{MD}}(\theta,m)$ is continuous (with respect to the $\|\cdot\|_\infty$-norm) in $m$ at $m = m_\theta$.*

LEMMA B.2.    *For all sequences of positive numbers $\delta_n = o(1)$,*

$$\sup_{\theta \in \Theta, \|m - m_\theta\|_{\mathcal{M}} \leq \delta_n} \|G_{n\mathrm{MD}}(\theta, m) - G_{\mathrm{MD}}(\theta, m)\|_2 = o_p(1).$$

LEMMA B.3.    *For all $(x, e)$, the ordinary partial derivative in $\theta$ of $G_{\mathrm{MD}}(\theta, m_\theta)(x, e)$, denoted $\Gamma_{1\mathrm{MD}}(\theta, m_\theta)(x, e)$, exists in a neighborhood of $\theta_o$, is continuous at $\theta = \theta_o$, and the matrix $\Gamma_{1\mathrm{MD}}(x, e) = \Gamma_{1\mathrm{MD}}(\theta_o, m_o)(x, e)$ is of full (column) rank for a set of positive $\mu$-measure $(x, e)$.*

LEMMA B.4.    *For $\mu$-all $(x, e)$, the pathwise derivative $\Gamma_{2\mathrm{MD}}(\theta, m_\theta)(x, e)$ of $G_{\mathrm{MD}}(\theta, m_\theta)(x, e)$ exists in all directions $m - m_\theta$ and satisfies the following:*

(i)    $\|G_{\mathrm{MD}}(\theta, m) - G_{\mathrm{MD}}(\theta, m_\theta) - \Gamma_{2\mathrm{MD}}(\theta, m_\theta)[m - m_\theta]\|_2 \leq c\|m - m_\theta\|_{\mathcal{M}}^2$

*for all $\theta$ with $\|\theta - \theta_o\| = o(1)$, all $m$ with $\|m - m_\theta\|_{\mathcal{M}} = o(1)$, some constant $c < \infty$;*

(ii)    $\|\Gamma_{2\mathrm{MD}}(\theta, m_\theta)[\widehat{m}_\theta - m_\theta] - \Gamma_{2\mathrm{MD}}(\theta_o, m_o)[\widehat{m} - m_o]\|_2$

$$\leq c\|\theta - \theta_o\| \times o_p(1) + O_p(n^{-1/2})$$

*for all $\theta$ with $\|\theta - \theta_o\| = o(1)$.*

LEMMA B.5.    *For all sequences of positive numbers $\{\delta_n\}$ with $\delta_n = o(1)$,*

$$\sup_{\|\theta - \theta_o\| \leq \delta_n, \|m - m_\theta\|_{\mathcal{M}} \leq \delta_n} \|G_{n\mathrm{MD}}(\theta, m) - G_{\mathrm{MD}}(\theta, m) - G_{n\mathrm{MD}}(\theta_o, m_o)\|_2$$

$$= o_p(n^{-1/2}).$$

LEMMA B.6.

$$\sqrt{n} \int \{G_{n\mathrm{MD}}(\theta_o, m_o) + \Gamma_{2\mathrm{MD}}(\theta_o, m_o)[\widehat{m} - m_o]\}(x, e)\Gamma_{1\mathrm{MD}}(x, e) \, d\mu(x, e)$$

$$\Longrightarrow N(0, V_{1\mathrm{MD}}).$$

**B.3. Proofs of Lemmas B.1–B.6.**

PROOF OF LEMMA B.1.    This follows from the representation

$$
\begin{aligned}
(22) \quad & G_{\mathrm{MD}}(\theta, m_\theta)(x, e) \\
& = E[[1(X \leq x) - F_X(x)]F_\varepsilon[\Lambda_o(\Lambda_\theta^{-1}(m_\theta(X) + e)) - m_o(X)]],
\end{aligned}
$$

and the smoothness of $F_\varepsilon$, $\Lambda_o$ and $\Lambda_\theta^{-1}$.    □

PROOF OF LEMMA B.2. Define the linearization

$$G_{n\text{MD}}^L(\theta, m)(x, e) = \widehat{F}_{X, \varepsilon(\theta, m)}(x, e) - F_X(x)\widehat{F}_{\varepsilon(\theta, m)}(e)$$
$$- \widehat{F}_X(x)F_{\varepsilon(\theta, m)}(e) + F_X(x)F_{\varepsilon(\theta, m)}(e).$$

By the triangle inequality, we have

$$\sup_{\theta \in \Theta, \|m - m_\theta\|_\mathcal{M} \leq \delta_n} \|G_{n\text{MD}}(\theta, m) - G_{\text{MD}}(\theta, m)\|_2$$

$$\leq \sup_{\theta \in \Theta, \|m - m_\theta\|_\mathcal{M} \leq \delta_n} \|G_{n\text{MD}}^L(\theta, m) - G_{\text{MD}}(\theta, m)\|_2$$

$$+ \sup_{\theta \in \Theta, \|m - m_\theta\|_\mathcal{M} \leq \delta_n} \|G_{n\text{MD}}(\theta, m) - G_{n\text{MD}}^L(\theta, m)\|_2.$$

We must show that both terms on the right-hand side are $o_p(1)$. Define the stochastic processes

$$\tau_{n\varepsilon}(\theta, m, e) = \widehat{F}_{\varepsilon(\theta, m)}(e) - F_{\varepsilon(\theta, m)}(e)$$

and

$$\tau_{nX\varepsilon}(\theta, m, x, e) = \widehat{F}_{X, \varepsilon(\theta, m)}(x, e) - F_{X, \varepsilon(\theta, m)}(x, e)$$

for each $\theta \in \Theta$, $m \in \mathcal{M}$, $x \in \mathbb{R}^k, e \in \mathbb{R}$. We claim that

$$(23) \qquad \sup_{\theta \in \Theta, \|m - m_\theta\|_\mathcal{M} \leq \delta_n, e \in \mathbb{R}} |\tau_{n\varepsilon}(\theta, m, e)| = o_p(1),$$

$$(24) \qquad \sup_{\theta \in \Theta, \|m - m_\theta\|_\mathcal{M} \leq \delta_n, x \in \mathbb{R}^k, e \in \mathbb{R}} |\tau_{nX\varepsilon}(\theta, m, x, e)| = o_p(1),$$

which implies that

$$\sup_{\theta \in \Theta, \|m - m_\theta\|_\mathcal{M} \leq \delta_n} \|G_{n\text{MD}}^L(\theta, m) - G_{\text{MD}}^L(\theta, m)\|_2$$

$$= \sup_{\theta \in \Theta, \|m - m_\theta\|_\mathcal{M} \leq \delta_n} \|(\widehat{F}_{X, \varepsilon(\theta, m)} - F_{X, \varepsilon(\theta, m)})$$

$$- F_X(\widehat{F}_{\varepsilon(\theta, m)} - F_{\varepsilon(\theta, m)}) - F_{\varepsilon(\theta, m)}(\widehat{F}_X - F_X)\|_2$$

$$\leq \Big[ \sup_{\theta \in \Theta, \|m - m_\theta\|_\mathcal{M} \leq \delta_n, e \in \mathbb{R}} |\tau_{nX\varepsilon}(\theta, m, e)|$$

$$+ \sup_{\theta \in \Theta, \|m - m_\theta\|_\mathcal{M} \leq \delta_n, x \in \mathbb{R}^k, e \in \mathbb{R}} |\tau_{n\varepsilon}(\theta, m, x, e)| + \sup_{x \in \mathbb{R}^k} |\widehat{F}_X(x) - F_X(x)| \Big]$$

$$= o_p(1).$$

Similarly, $\sup_{\theta \in \Theta, \|m - m_\theta\|_\mathcal{M} \leq \delta_n} \|G_{n\text{MD}}(\theta, m) - G_{n\text{MD}}^L(\theta, m)\|_2 = o_p(1)$. The proof of (23) and (24) is based on Theorem 3 in CLV (2003). We omit the details because it is similar to our proof of Lemma B.5. □

PROOF OF LEMMA B.3. Below, we calculate $\Gamma_{1\mathrm{MD}}(x,e) = \Gamma_{1\mathrm{MD}}(\theta_o, m_o) \times (x,e)$. In a similar way $\Gamma_{1\mathrm{MD}}(\theta, m_\theta)(x,e)$ can be obtained. First, we have

$$\frac{\partial}{\partial \theta} F_{\varepsilon(\theta, m_\theta)}(e) \Big\downarrow_{\theta=\theta_o}$$

$$= E \frac{\partial}{\partial \theta} F_\varepsilon[\Lambda_o(\Lambda_\theta^{-1}(m_\theta(X)+e)) - m_o(X)] \downarrow_{\theta=\theta_o}$$

$$= f_\varepsilon(e) E \frac{\partial}{\partial \theta} \Lambda_o(\Lambda_\theta^{-1}(m_\theta(X)+e)) \Big\downarrow_{\theta=\theta_o}$$

$$= f_\varepsilon(e) E \Lambda_o'(\Lambda_o^{-1}(m_o(X)+e)) \frac{\partial}{\partial \theta} (\Lambda_\theta^{-1}(m_\theta(X)+e)) \Big\downarrow_{\theta=\theta_o}$$

$$= f_\varepsilon(e) E \Lambda_o'(\Lambda_o^{-1}(m_o(X)+e)) \left[ \frac{\lambda_o(\Lambda_o^{-1}(m_o(X)+e))}{\Lambda_o'(\Lambda_o^{-1}(m_o(X)+e))} \right.$$

$$\left. + \frac{1}{\Lambda_o'(\Lambda_o^{-1}(m_o(X)+e))} \dot{m}_o(X) \right]$$

$$= f_\varepsilon(e) E[\lambda_o(\Lambda_o^{-1}(m_o(X)+e)) + m_o(X)]$$

by the chain rule. Similarly,

$$\frac{\partial}{\partial \theta} F_{X,\varepsilon(\theta,m_\theta)}(x,e) \Big\downarrow_{\theta=\theta_o}$$

$$= f_\varepsilon(e) E[1(X \leq x)\{\lambda_o(\Lambda_o^{-1}(m_o(X)+e)) + \dot{m}_o(X)\}].$$

Therefore,

$$(25) \quad \begin{aligned} \Gamma_{1\mathrm{MD}}(x,e) &= \Gamma_{1\mathrm{MD}}(\theta_o, m_o)(x,e) = \frac{\partial G_{\mathrm{MD}}(\theta, m_\theta)}{\partial \theta}(x,e) \Big\downarrow_{\theta=\theta_o} \\ &= \frac{\partial}{\partial \theta} F_{X,\varepsilon(\theta,m_\theta)}(x,e) - F_X(x) \frac{\partial}{\partial \theta} F_{\varepsilon(\theta,m_\theta)}(e) \\ &= f_\varepsilon(e) E[(1(X \leq x) - F_X(x)) \\ &\qquad \times (\lambda_o(\Lambda_o^{-1}(m_o(X)+e)) + \dot{m}_o(X))]. \quad \square \end{aligned}$$

PROOF OF LEMMA B.4. By the law of iterated expectation and partial differentiation, we obtain that

$$[\Gamma_{2\mathrm{MD}}(\theta_o, m_o)(m-m_o)](x,e)$$

$$= \frac{\partial G_{\mathrm{MD}}(\theta_o, m_o + t(m-m_o))}{\partial t}(x,e) \Big\downarrow_{t=0}$$

$$= f_\varepsilon(e) E[(1(X \leq x) - F_X(x))(m(X) - m_o(X))].$$

Similarly, the formula of $[\Gamma_{2\mathrm{MD}}(\theta, m_\theta)(m - m_\theta)](x, e)$ is given by

$$[\Gamma_{2\mathrm{MD}}(\theta, m_\theta)(m - m_\theta)](x, e)$$
$$= \lim_{\tau \to 0} \frac{1}{\tau} E[\{1(X \leq x) - F_X(x)\} f_\varepsilon[\Lambda_o\{\Lambda_\theta^{-1}(m_\theta(X) + e)\} - m_o(X)]$$
$$\times [\Lambda_o\{\Lambda_\theta^{-1}(m_\theta(X) + \tau(m - m_\theta)(X) + e)\}$$
$$- \Lambda_o\{\Lambda_\theta^{-1}(m_\theta(X) + e)\}]].$$

The two inequalities in the statement of Lemma B.4 now follow easily, using the consistency of $\widehat{m}_\theta$ and the fact that $\sup_x |(\widehat{m}_\theta - m_\theta)(x) - (\widehat{m}_o - m_o)(x)| = o_p(1)\|\theta - \theta_o\| + O_p(n^{-1/2})$. □

PROOF OF LEMMA B.5. Define the stochastic processes

$$\nu_{n\varepsilon}(\theta, m, e) = \sqrt{n}[\widehat{F}_{\varepsilon(\theta,m)}(e) - F_{\varepsilon(\theta,m)}(e)]$$

and

$$\nu_{nX\varepsilon}(\theta, m, x, e) = \sqrt{n}[\widehat{F}_{X, \ \varepsilon(\theta,m)}(x, e) - F_{X,\varepsilon(\theta,m)}(x, e)]$$

for each $\theta \colon \|\theta - \theta_o\| \leq \delta_n$ and $m \colon \|m - m_\theta\|_{\mathcal{M}} \leq \delta_n$, $x \in \mathbb{R}^k$, $e \in \mathbb{R}$. We claim that

$$(26) \qquad \sup_{\|\theta - \theta_o\| \leq \delta_n, \|m - m_\theta\|_{\mathcal{M}} \leq \delta_n, e \in \mathbb{R}} |\nu_{n\varepsilon}(\theta, m, e)| = o_p(1),$$

$$(27) \qquad \sup_{\|\theta - \theta_o\| \leq \delta_n, \|m - m_\theta\|_{\mathcal{M}} \leq \delta_n, x \in \mathbb{R}^d, e \in \mathbb{R}} |\nu_{nX\varepsilon}(\theta, m, x, e)| = o_p(1).$$

The proof of these results is based on Theorem 3 in CLV (2003). We have to show that their condition (3.2) is satisfied, which requires in our case [with $g(Z, \theta, m) = 1(\varepsilon(\theta, m) \leq e) - E1(\varepsilon(\theta, m) \leq e)$ and $g(Z, \theta, m) = 1(X \leq x)1(\varepsilon(\theta, m) \leq e) - E1(X \leq x)1(\varepsilon(\theta, m) \leq e)$] that

$$\left( E\left[ \sup_{(\theta', m') \colon \|\theta' - \theta\| < \delta, \|m' - m\|_{\mathcal{M}} < \delta} |g(Z, \theta', m') - g(Z, \theta, m)|^r \right] \right)^{1/r} \leq K\delta^s$$

for all $(\theta, m) \in \Theta \times \mathcal{M}$, all small positive value $\delta = o(1)$, and for some constants $s \in (0, 1]$, $K > 0$, and that the bound holds for $\mu$-almost all $(x, e)$. We have

$$|g(Z, \theta', m') - g(Z, \theta, m)| \leq |1(\varepsilon(\theta, m) \leq e) - 1(\varepsilon(\theta', m') \leq e)|$$
$$+ |E1(\varepsilon(\theta, m) \leq e) - E1(\varepsilon(\theta', m') \leq e)|$$

and

$$|1(\varepsilon(\theta, m) \leq e) - 1(\varepsilon(\theta', m') \leq e)|$$
$$= |1(\Lambda_\theta(Y) - m(X) \leq e) - 1(\Lambda_{\theta'}(Y) - m'(X) \leq e)|$$
$$\leq |1(\Lambda_\theta(Y) - m(X) \leq e) - 1(\Lambda_\theta(Y) - m'(X) \leq e)|$$
$$+ |1(\Lambda_\theta(Y) - m'(X) \leq e) - 1(\Lambda_{\theta'}(Y) - m'(X) \leq e)|.$$

For all $m' \in \mathcal{M}$ with $\|m' - m\|_{\mathcal{M}} \leq \delta \leq 1$, we have for all $Y$, $X$, $e$

$$\sup_{\|m'-m\|_{\mathcal{M}} \leq \delta} |1(m'(X) \geq \Lambda_\theta(Y) - e) - 1(m(X) \geq \Lambda_\theta(Y) - e)|$$

$$\leq 1(m(X) + \delta \geq \Lambda_\theta(Y) - e) - 1(m(X) - \delta \geq \Lambda_\theta(Y) - e).$$

The preceding term is either one or zero and its expectation is the probability that $m(X) + \delta \geq \Lambda_\theta(Y) - e \geq m(X) - \delta$, which is the probability that $e + \delta \geq \Lambda_\theta(Y) - m(X) \geq e - \delta$, which is

$$F_{\varepsilon(\theta,m)}(e + \delta) - F_{\varepsilon(\theta,m)}(e - \delta)$$

$$= EF_\varepsilon[\Lambda_o(\Lambda_\theta^{-1}(m(X) + e + \delta)) - m_o(X)]$$

$$- EF_\varepsilon[\Lambda_o(\Lambda_\theta^{-1}(m(X) + e - \delta)) - m_o(X)].$$

We then apply the smoothness conditions on $F_\varepsilon$, $\Lambda_o$ and $\Lambda_\theta^{-1}$ to bound the right-hand side by $K\delta$ for small enough $\delta$ and constant $K < \infty$.

Next, by the Mean Value Theorem, we have $\Lambda_\theta(Y) - \Lambda_{\theta'}(Y) = \lambda_{\theta^*}(Y) \times (\theta - \theta')$, where $\theta^*$ is an intermediate value between $\theta$ and $\theta'$. For all $\alpha > 0$, by the Bonferroni and Markov inequalities,

$$\Pr\left[\max_{1 \leq i \leq n} \sup_{\|\theta - \theta'\| \leq \delta} |\lambda_{\theta'}(Y_i)| > c \times n^\alpha\right]$$

$$\leq n \times \Pr\left[\sup_{\|\theta - \theta'\| \leq \delta} |\lambda_{\theta'}(Y)| > c \times n^\alpha\right]$$

$$\leq n \times \frac{E[\sup_{\|\theta - \theta'\| \leq \delta} |\lambda_{\theta'}(Y)|^k]}{c^k n^{k\alpha}} = o(1),$$

provided $k > \alpha^{-1}$.

Therefore, we can safely assume that there is some upper bound $c$ such that $\sup_{\|\theta - \theta'\| \leq \delta} |\Lambda_\theta(Y) - \Lambda_{\theta'}(Y)| \leq c \times \delta$. Therefore, on this set,

$$\sup_{\|\theta' - \theta\| \leq \delta} |1(\Lambda_\theta(Y) - m'(X) \leq e) - 1(\Lambda_{\theta'}(Y) - m'(X) \leq e)|$$

$$\leq 1(\Lambda_\theta(Y) + c\delta - m'(X) \leq e) - 1(\Lambda_\theta(Y) - c\delta - m'(X) \leq e)|,$$

which has probability bounded by $K\delta$ for some $K > 0$.

Therefore, condition (3.2) of Theorem 3 in CLV (2003) is satisfied with $r = 2$ and $s = 1/2$, and condition (3.3) of Theorem 3 is satisfied by the condition on the covering number of the class $\mathcal{M}$, stated in Assumption B.8.

$\square$

PROOF OF LEMMA B.6. We show below that

$$[\Gamma_{2\text{MD}}(\theta_o, m_o)(\widehat{m} - m_o)](x, e)$$

(28)
$$= f_\varepsilon(e)\sqrt{n}\int [(1(X \le x) - F_X(x))(\widehat{m}(X) - m_o(X))]f_X(X)\,dX$$

$$= f_\varepsilon(e)\frac{1}{\sqrt{n}}\sum_{i=1}^n (1(X_i \le x) - F_X(x))$$

$$\times f_X(X_i)\sum_{\alpha=1}^d v_{o1\alpha}(X_{\alpha i}, \varepsilon_i) + o_p(1).$$

Therefore,

$$[L_n(x, e) + [\Gamma_{2\text{MD}}(\theta_o, m_o)(\widehat{m} - m_o)](x, e)] = \frac{1}{\sqrt{n}}\sum_{i=1}^n U_i(x, e) + o_p(1),$$

where

$$U_i(x, e) = [1(X_i \le x)1(\varepsilon_i \le e) - F_{X,\varepsilon}(x, e)]$$

$$- F_X(x)[1(\varepsilon_i \le e) - F_\varepsilon(e)]$$

$$- F_\varepsilon(e)[1(X_i \le x) - F_X(x)]$$

$$+ f_X(X_i)\sum_{\alpha=1}^d v_{o1\alpha}(X_{\alpha i}, \varepsilon_i)f_\varepsilon(e)(1(X_i \le x) - F_X(x)),$$

and where $E[U_i(x, e)] = 0$ for all $x, e$. Because $F_{X,\varepsilon}(x, e) = F_X(x)F_\varepsilon(e)$, we have

$$U_i(x, e) = [1(X_i \le x) - F_X(x)][1(\varepsilon_i \le e) - F_\varepsilon(e)]$$

$$+ f_X(X_i)\sum_{\alpha=1}^d v_{o1\alpha}(X_{\alpha i}, \varepsilon_i)f_\varepsilon(e)(1(X_i \le x) - F_X(x)).$$

Now integrating $U_i(x, e)$ with respect to $\Gamma_{1\text{MD}}(x, e)\,d\mu(x, e)$ gives the answer.

*Proof of* (28): Write

$$\widehat{m}(X) - m_o(X)$$

$$= \frac{1}{nh}\sum_{i=1}^n\sum_{\alpha=1}^d K_1\left(\frac{X_\alpha - X_{\alpha i}}{h}\right)v_{o1\alpha}(X_{\alpha i}, \varepsilon_i)$$

$$+ \frac{1}{n}\sum_{i=1}^n v_{o2}(\varepsilon_i) + o_p(n^{-1/2}).$$

Then, provided $nh^{2q_1} \to 0$,

$$\sqrt{n} \int [(1(X \le x) - F_X(x))(\hat{m}(X) - m_o(X))] f_X(X) \, dX$$

$$= \frac{1}{\sqrt{n}} \sum_{i=1}^n \sum_{\alpha=1}^d v_{o1\alpha}(X_{\alpha i}, \varepsilon_i)$$

$$\times \int \left[ (1(X \le x) - F_X(x)) \frac{1}{h} K_1 \left( \frac{X_\alpha - X_{\alpha i}}{h} \right) \right] f_X(X) \, dX$$

$$+ \frac{1}{\sqrt{n}} \sum_{i=1}^n v_{o2}(\varepsilon_i) \int [(1(X \le x) - F_X(x))] f_X(X) \, dX + o_p(1)$$

$$= \frac{1}{\sqrt{n}} \sum_{i=1}^n \sum_{\alpha=1}^d v_{o1\alpha}(X_{\alpha i}, \varepsilon_i) \int [(1(X_i + uh \le x) - F_X(x)) K_1(u_\alpha)]$$

$$\times f_X(X_i + uh) \, du + o_p(1)$$

$$= \frac{1}{\sqrt{n}} \sum_{i=1}^n \sum_{\alpha=1}^d v_{o1\alpha}(X_{\alpha i}, \varepsilon_i)(1(X_i \le x) - F_X(x)) f_X(X_i) + o_p(1).$$

We also have to substitute $\frac{\partial m_\theta}{\partial \theta}(x) \downarrow_{\theta=\theta_o}$ into the formula for $\Gamma_{1MD}$. $\quad \square$

**Acknowledgments.** We thank Enno Mammen and two anonymous referees for helpful discussion.

O. Linton
Department of Economics
London School of Economics
Houghton Street, London WC2A 2AE
United Kingdom
E-mail: o.linton@lse.ac.uk

S. Sperlich
Institut für Statistik
  und Ökonometrie
Georg-August Universität
Platz der Göttinger Sieben 5
37073 Göttingen
Germany
E-mail: stefan.sperlich@wiwi.uni-goettingen.de

I. Van Keilegom
Institut de Statistique
Université catholique de Louvain
Voie du Roman Pays 20
B 1348 Louvain-la-Neuve
Belgium
E-mail: vankeilegom@stat.ucl.ac.be